\pdfoutput=1
\RequirePackage{ifpdf}
\ifpdf 
\documentclass[pdftex]{sigma}
\else
\documentclass{sigma}
\fi

\usepackage{amscd,array,enumerate}

\newcommand{\onto}{\to\!\!\!\!\!\to}

\numberwithin{equation}{section}
\newtheorem{Theorem}{Theorem}[section]
\newtheorem{Conjecture}[Theorem]{Conjecture}
\newtheorem{Corollary}[Theorem]{Corollary}

\newtheorem{Lemma}[Theorem]{Lemma}
\newtheorem{Proposition}[Theorem]{Proposition}
{ \theoremstyle{definition}
\newtheorem{Definition}[Theorem]{Definition}
\newtheorem{Example}[Theorem]{Example}
\newtheorem{Remark}[Theorem]{Remark}
}

\begin{document}

\newcommand{\arXivNumber}{1401.4338}

\allowdisplaybreaks

\renewcommand{\thefootnote}{$\star$}

\renewcommand{\PaperNumber}{024}

\FirstPageHeading

\ShortArticleName{The Feigin Tetrahedron}

\ArticleName{The Feigin Tetrahedron\footnote{This paper is a~contribution to the Special Issue on New Directions in Lie
Theory. The full collection is available at \href{http://www.emis.de/journals/SIGMA/LieTheory2014.html}{http://www.emis.de/journals/SIGMA/LieTheory2014.html}}}

\Author{Dylan RUPEL}

\AuthorNameForHeading{D.~Rupel}

\Address{Department of Mathematics, Northeastern University, Boston, MA 02115, USA}
\Email{\href{mailto:d.rupel@neu.edu}{d.rupel@neu.edu}}
\URLaddress{\url{http://www.northeastern.edu/drupel/}}

\ArticleDates{Received September 11, 2014, in f\/inal form March 03, 2015; Published online March 19, 2015}

\Abstract{The f\/irst goal of this note is to extend the well-known Feigin homomorphisms taking quantum groups to quantum
polynomial algebras.
More precisely, we def\/ine genera\-li\-zed Feigin homomorphisms from a~quantum shuf\/f\/le algebra to quantum polynomial algebras
which extend the classical Feigin homomorphisms along the embedding of the quantum group into said quantum shuf\/f\/le
algebra.
In a~recent work of Berenstein and the author, analogous extensions of Feigin homomorphisms from the dual Hall--Ringel
algebra of a~valued quiver to quantum polynomial algebras were def\/ined.
To relate these constructions, we establish a~homomorphism, dubbed the \emph{quantum shuffle character}, from the dual
Hall--Ringel algebra to the quantum shuf\/f\/le algebra which relates the generalized Feigin homomorphisms.
These constructions can be compactly described by a~commuting tetrahedron of maps beginning with the quantum group and
terminating in a~quantum polynomial algebra.
The second goal in this project is to better understand the dual canonical basis conjecture for skew-symmetrizable
quantum cluster algebras.
In the symmetrizable types it is known that dual canonical basis elements need not have positive multiplicative
structure constants, while this is still suspected to hold for skew-symmetrizable quantum cluster algebras.
We propose an alternate conjecture for the symmetrizable types: the cluster monomials should correspond to irreducible
characters of a~KLR algebra.
Indeed, the main conjecture of this note would establish this ``KLR conjecture'' for acyclic skew-symmetrizable quantum
cluster algebras: that is, we conjecture that the images of rigid representations under the quantum shuf\/f\/le character
give irreducible characters for KLR algebras.
We sketch a~proof in the symmetric case giving an alternative to the proof of Kimura--Qin that all non-initial cluster
variables in an acyclic skew-symmetric quantum cluster algebra are contained in the dual canonical basis.
With these results in mind we interpret the cluster mutations directly in terms of the representation theory of the KLR
algebra.}


\Classification{13F60; 16G20; 17B37; 20G42}

\renewcommand{\thefootnote}{\arabic{footnote}} 
\setcounter{footnote}{0}

\section{Introduction}

In the struggle to understand the quantum groups and their duals, the quantized enveloping algebras, one searches for
concrete realizations of these algebras, either by generators and relations or via embeddings, or more generally simply
homomorphisms, into ``nicer'' algebras.
This note grew out of an attempt to understand three well-known and well-studied examples, namely the embedding into
a~(dual) Hall--Ringel algebra, the embedding into a~quantum shuf\/f\/le algebra, and the Feigin homomorphisms to quantum
polynomial algebras.

The embeddings of the (positive part of) quantized Kac--Moody algebras into Hall--Ringel algebras originated with the
incredible insights of Ringel, for the f\/inite-types in~\cite{ringel2} and fully realized in~\cite{ringel3}.
This is merely the tip of a~fantastic geometric iceberg discovered by Lusztig~\cite{lusztig} that has lead to a~deep
understanding of certain ``canonical'' bases in the quantized Kac--Moody algebras.

A careful study of multiplicative properties of these bases helped motivate Fomin and Zele\-vinsky in the def\/inition of
cluster algebras.
In the hopes of gaining a~combinatorial grasp of the structure of the canonical basis they provided the def\/inition of
a~recursively, combinatorially def\/ined algebra with a~deep conjecture in mind: certain ``cluster monomials'' appearing
in this construction should be identif\/iable with (dual) canonical basis vectors.
This motivating conjecture has been established so far in a~very restricted set of cases by Lampe~\cite{lampe1,lampe2}
and Kimura--Qin~\cite{kimura-qin}.
One aim of this note is to shed some additional light on this ``dual canonical basis conjecture''.

The quantized Kac--Moody algebras are most compactly described by generators and relations or better as a~certain
quotient of a~free algebra with an unconventional (dare I say twisted) multiplication def\/ined on its tensor square.
Dualizing this construction leads to an embedding of the quantum group into the quantum shuf\/f\/le algebra and thus, via an
isomorphism, an embedding of the quantized Kac--Moody algebra.
This is the tip of yet another (not entirely unrelated) iceberg, this one
categorical~\cite{khovanov-lauda1,khovanov-lauda2,rouquier1,rouquier2}, which provides another description of these
positively amazing bases of the quantum group and quantized Kac--Moody algebra.

Our third major tool for studying quantized Kac--Moody algebras are certain algebra homomorphisms to quantum polynomial
rings suggested by B.~Feigin during a~talk at RIMS in 1992 (thus we refer to these as Feigin homomorphisms).
There are actually inf\/initely many Feigin homomorphisms acting on a~given quantized Kac--Moody algebra, one
corresponding to each f\/inite sequence of simple roots.
For f\/inite-types, we get an embedding whenever the sequence corresponds to a~reduced word for the longest element of the
corresponding Weyl group.
In general a~Feigin homomorphism has a~very large kernel, a~feature which we hope will shed some light on the
relationship of KLR characters to quantum cluster algebras.
The Feigin homomorphisms, as Feigin predicted, turn out to be the essential tool~\cite{berenstein,iohara-malikov, joseph} for studying the skew-f\/ield of fractions of the quantized Kac--Moody algebras.

The relationship between the f\/irst and third approach was studied by Berenstein and the author
in~\cite{berenstein-rupel} where Feigin homomorphisms were extended to the dual Hall--Ringel algebra for representations  
of an acyclic valued quiver.
Applications to quantum cluster algebras were a~central result in~\cite{berenstein-rupel} where the images of rigid  
objects are identif\/ied with quantum cluster characters describing non-initial cluster variables of a~corresponding
acyclic quantum cluster algebra~\cite{qin,rupel1,rupel2}.
This also provided the essential tool for understanding a~certain ``twist automorphism'' of the quantum cluster algebra,
it would be interesting to use the constructions/conjectures of this note to interpret the twist in terms of the
representation theory of KLR algebras.

The main result of the current work is to complete the above tetrahedron.
We begin in Section~\ref{sec:shuffle} by extending Feigin homomorphisms to the quantum shuf\/f\/le algebra.
Section~\ref{sec:shuffle_character} introduces a~homomorphism, the \emph{quantum shuffle character}, from the dual
Hall--Ringel algebra to the quantum shuf\/f\/le algebra making all of the relevant triangles commute.

A secondary goal of this project is to understand the dual canonical basis conjecture for quantum cluster algebras.
The generalized Feign homomorphism of~\cite{berenstein-rupel} produces non-initial cluster monomials (of a~certain
acyclic quantum cluster algebra) as the images of rigid objects.
Tracing the other path through the tetrahedron gives elements of the quantum shuf\/f\/le algebra which, according to the
dual canonical basis conjecture, should be expected to identify (in symmetric types) with dual canonical basis elements,
i.e.~with irreducible characters of Khovanov--Lauda--Rouquier (quiver Hecke)
algebras~\cite{khovanov-lauda1,khovanov-lauda2,rouquier1,rouquier2} living in the quantum shuf\/f\/le algebra.
In Section~\ref{sec:KLR} we recall the necessary background material to precisely formulate this conjecture as well as
sketch a~proof for symmetric types giving an alternative to the proof by Kimura and Qin using graded quiver varieties.
Assuming the conjecture, we formulate a~mutation operation for certain clusters def\/ined using the representation theory
of KLR algebras.
We f\/inish by presenting several conjectures about how this should be done in general.

The results of this work are undoubtedly closely related to several other works in this area though the precise
connections are still a~little bit of a~mystery, perhaps the other constructions provide the categorif\/ied/geometric
perspective of our constructions.
We mention here only a~few main examples that should be investigated further.
The geometric construction of KLR algebras by Varagnolo--Vasserot~\cite{varagnolo-vasserot} and Webster~\cite{webster}
motivate our approach to the dual canonical basis conjecture in acyclic skew-symmetric types, more precisely they have
given a~geometric construction of certain faithful modules over KLR algebra using f\/lag varieties of quiver
representations.
The other works we mention are mainly focused on understanding the precise relationship of the representation theory of
KLR algebras to an alternative categorif\/ication of quantum groups using quantum af\/f\/ine
algebras~\cite{kang-kashiwara-kim}.
A~central ingredient in establishing a~direct link was the paper~\cite{hernandez-leclerc2} where the authors f\/ind an  
isomorphism between between a~deformed Grothendieck ring associated to the quantum af\/f\/ine algebra and the derived Hall
algebra of a~quiver~\cite{toen}.
Tracing through these isomorphisms one should f\/ind a~direct relationship between the derived Hall algebra and the KLR
algebra, perhaps a~Feigin-type homomorphism is lurking in the background.

The structure of the paper is as follows: Section~\ref{sec:combinatorics} sets up certain combinatorial and notational
conventions that we have found most useful in these investigations.
Section~\ref{sec:feigin} recalls the quantized Kac--Moody algebra and background on classical Feigin homomorphisms.
In Section~\ref{sec:qca} we introduce quantum cluster algebras.
Section~\ref{sec:valued quivers} def\/ines valued quivers and their Hall--Ringel algebras and recalls the results
of~\cite{qin,rupel1,rupel2} and~\cite{berenstein-rupel} relating these to quantum cluster algebras.
Sections~\ref{sec:shuffle}--\ref{sec:KLR} have been discussed above.
Section~\ref{sec:generalizations} discusses ideas for generalizations and proposes several conjectures on irreducible KLR modules.

\section{Combinatorial conventions and notations}\label{sec:combinatorics}

Fix an indeterminate~$q$.
We def\/ine the~$q$-numbers,~$q$-factorials, and~$q$-binomials~by
\begin{gather*}
(n)=(n)_q=1+q+\cdots+q^{n-1},
\\
(n)!=(n)\cdot(n-1)\cdots(2)\cdot(1),
\qquad
{n\choose k}_q=\frac{(n)!}{(k)!\cdot(n-k)!}.
\end{gather*}
In def\/ining our generalized Feigin homomorphism $\overline\Psi_\mathbf{i}$ we will make particular use of the
bar-invariant versions of the quantum numbers, for convenience we introduce $v=\sqrt{q}$:
\begin{gather*}
[n]=[n]_q=v^{-n+1}+v^{-n+3}+\cdots+v^{n-3}+v^{n-1},
\\
[n]!=[n]\cdot[n-1]\cdots[2]\cdot[1],
\qquad
{n\brack k}_q=\frac{[n]!}{[k]!\cdot[n-k]!}.
\end{gather*}
The~$q$-binomials satisfy several useful analogues of classical binomial identities.
We collect these and sketch their proofs in the following.
\begin{Lemma}\label{lem:binomial_identities}\qquad
\begin{enumerate}\itemsep=0pt
\item[$1.$] The~$q$-binomial coefficients satisfy the following identities:
\begin{enumerate}[\upshape (a)]\itemsep=0pt
\item Pascal identities: ${n\choose k}_q={n-1\choose k-1}_q+q^k{n-1\choose k}_q=q^{n-k}{n-1\choose k-1}_q+{n-1\choose k}_q$;
\item row identity: for $n>0$, $\sum\limits_{k=0}^n (-1)^k q^{-nk+\frac{1}{2} k(k+1)}{n\choose k}_q=0$;
\item subspace identity: ${m+n\choose k}_q=\sum\limits_{r+s=k}q^{r(n-s)}{m\choose r}_q{n\choose s}_q$.
\end{enumerate}
\item[$2.$] The bar-invariant~$q$-binomial coefficients satisfy the following identities:
\begin{enumerate}[\upshape (a)]\itemsep=0pt
\item Pascal identities: ${n\brack k}_q=v^{n-k}{n-1\brack k-1}_q+v^{-k}{n-1\brack k}_q=v^{k-n}{n-1\brack k-1}_q+v^k{n-1\brack k}_q$;
\item row identity: for $n>0$, $\sum\limits_{k=0}^n (-1)^k v^{k-nk} {n\brack k}_q=0$;
\item subspace identity: ${m+n\brack k}_q=\sum\limits_{r+s=k}v^{r(n-s)-s(m-r)}{m\brack r}_q{n\brack s}_q$.
\end{enumerate}
\end{enumerate}
\end{Lemma}

\begin{proof}
We prove the identities in (1).
The Pascal identities in (a) follow from the def\/inition and the equalities
$(n)_q=(k)_q+q^k(n-k)_q=q^{n-k}(k)_q+(n-k)_q$.
For (b) we simply apply the second Pascal identity to get a~telescopic summation:
\begin{gather*}
\sum\limits_{k=0}^n (-1)^k q^{-nk+\frac{1}{2} k(k+1)}{n\choose k}_q
\\
\qquad
=\sum\limits_{k=1}^n (-1)^k
q^{-n(k-1)+\frac{1}{2}(k-1)k}{n-1\choose k-1}_q+\sum\limits_{k=0}^{n-1} (-1)^k q^{-nk+\frac{1}{2} k(k+1)}{n-1\choose
k}_q=0.
\end{gather*}
To see (c) we consider a~f\/inite f\/ield $\mathbb{F}$ with~$q$ elements.
Write $\operatorname{Gr}_k(\mathbb{F}^n)$ for the set of~$k$-dimensional subspaces of $\mathbb{F}^n$.
The number of points in $\operatorname{Gr}_k(\mathbb{F}^n)$ is given by ${n\choose k}_q$.
Fix a~decomposition $\mathbb{F}^{m+n}=\mathbb{F}^m\oplus\mathbb{F}^n$ so that any subspace of $\mathbb{F}^{m+n}$ can be
decomposed as a~direct sum of a~subspace of $\mathbb{F}^m$ and a~subspace of $\mathbb{F}^n$.
This gives rise to a~surjective map $\operatorname{Gr}_k(\mathbb{F}^{m+n})\onto \bigsqcup\limits_{r+s=k}\operatorname{Gr}_r(\mathbb{F}^m)\times\operatorname{Gr}_s(\mathbb{F}^n)$ with f\/iber over any point
of $\operatorname{Gr}_r(\mathbb{F}^m)\times\operatorname{Gr}_s(\mathbb{F}^n)$ an af\/f\/ine space of dimension $r(n-s)$.
Now counting points completes the proof when~$q$ is a~power of a~prime.
But this is an equality of polynomials for inf\/initely many values and thus the polynomials must be equal for every~$q$.

The identities in (2) easily follow from those in (1) using the relation ${n\brack k}_q=v^{-k(n-k)}{n\choose k}_q$.
\end{proof}

We choose a~notation for symmetric groups which is most convenient for describing the structure constants of the
multiplication in the quantum shuf\/f\/le algebra.
For each positive integer~$t$ let $\Sigma_t$ denote the symmetric group on~$t$ letters, thought of as the set of all
orderings $\sigma=(\sigma_1,\ldots,\sigma_t)$ of the set $(1,\ldots,t)$.
We write $\sigma\cdot\tau=(\sigma_{\tau_1},\ldots,\sigma_{\tau_t})$ for the multiplication in $\Sigma_t$.
Let $\sigma^{-1}_k$ denote the position of~$k$ in~$\sigma$, i.e the number~$i$ so that $\sigma_i=k$.

For positive integers~$r$ and~$s$ def\/ine the \emph{shuffle subgroup} $\Sigma_{r,s}\subset\Sigma_{r+s}$ consisting of all
shuf\/f\/les of the sequences $(1,\ldots,r)$ and $(r+1,\ldots,r+s)$, that is
\begin{gather*}
\Sigma_{r,s}=\big\{\sigma\in\Sigma_{r+s}: \sigma^{-1}_1<\cdots<\sigma^{-1}_r, \sigma^{-1}_{r+1}<\cdots<\sigma^{-1}_{r+s}\big\}.
\end{gather*}
By convention we take $\Sigma_{0,0}=\Sigma_0=\{1\}$ to be the trivial group.

\section{Quantum groups and Feigin homomorphisms}\label{sec:feigin}

In this section we introduce the main objects underlying the present investigations: the quantum groups associated to
a~symmetrizable Cartan matrix and the Feigin homomorphisms mapping them to quantum polynomial algebras.

Fix an index set~$I$.
Let $A=(a_{ij})$ be an $I\times I$ Cartan matrix with symmetrizing matrix $D=\operatorname{diag}(\mathbf{d})$, where
$\mathbf{d}=(d_i: i\in I)$ is an~$I$-tuple of positive integers.
More explicitly, the entries of~$A$ have the following properties:
\begin{itemize}\itemsep=0pt
\item $a_{ii}=2$ for all~$i$;
\item $a_{ij}\le0$ for $i\ne j$;
\item $d_ia_{ij}=d_ja_{ji}$ for all $i,j\in I$.
\end{itemize}
Let $\mathcal{Q}$ denote the root lattice of a~Kac--Moody root system~$\Phi$ associated to~$A$.
Denoting by $\{\alpha_i\}_{i\in I}$ the simple roots of~$\Phi$, $\mathcal{Q}$ is the free abelian 
group with basis $\{\alpha_i\}$.
The pair $(A,\mathbf{d})$ determines a~symmetric bilinear form $(\cdot,\cdot):\mathcal{Q}\times\mathcal{Q}\to\mathbb{Z}$
given on generators of $\mathcal{Q}$ by $(\alpha_i,\alpha_j)=d_ia_{ij}$ for $i,j\in I$.
Let~$W$ denote the Weyl group associated to~$\Phi$ with generators the simple ref\/lections $s_i$ ($i\in I$).
Write $\alpha_i^\vee=\alpha_i/d_i$.
The action of a~simple generator $s_i$ on a~root~$\alpha$ is given~by
\begin{gather*}
s_i(\alpha)=\alpha-(\alpha_i^\vee,\alpha)\alpha_i.
\end{gather*}
Notice that $(\cdot,\cdot)$ is equivariant with respect to the action of~$W$, i.e.~$(\alpha,\beta)=(s_i\alpha,s_i\beta)$
for all $\alpha,\beta\in\mathcal{Q}$ and $i\in I$.

Fix an indeterminate~$q$ and write $v=\sqrt{q}$.
Def\/ine $\mathcal{A}_v=\mathbb{Z}[v^{\pm1}]$.
It will be convenient for $i\in I$ to abbreviate $(n)_i=(n)_{q^{d_i}}$ and def\/ine $(n)^!_i$, ${n\choose k}_i$ similarly.
Similar abbreviations will be used for the bar-invariant versions.
Associated to the Cartan matrix~$A$ with choice of symmetrizers $\mathbf{d}$ we have the (integral form of the)
quantized coordinate ring $\mathcal{A}_v[N]$ of the standard upper unipotent subgroup~$N$ of the Kac--Moody group~$G$.
The $\mathcal{Q}$-graded algebra $\mathcal{A}_v[N]$ is generated by formal variables $x_i^{[r]}$ ($i\in I$,
$r\in\mathbb{Z}_{\ge0}$) with $|x_i^{[r]}|=-r\alpha_i$ subject to the quantum Serre relations:
\begin{gather*}
\sum\limits_{r=0}^{1-a_{ij}} (-1)^r x_i^{[r]}x_jx_i^{[1-a_{ij}-r]}=0
\qquad
\text{for any} \ \ i\ne j\in I,
\end{gather*}
where by def\/inition $[r]_i^!x_i^{[r]}=x_i^r$.
The algebra $\mathcal{A}_v[N]$ is actually a~twisted bialgebra where  
each generator is primitive, i.e.~the comultiplication is given
by $\Delta(x_i)=x_i\otimes1+1\otimes x_i$ for each $i\in I$, and where the multiplication on  
$\mathcal{A}_v[N]\otimes\mathcal{A}_v[N]$ is twisted by the grading, i.e.~$(x\otimes y)(x'\otimes y')=v^{(|x'|,|y|)}\cdot xx'\otimes yy'$.

The graded dual of $\mathcal{A}_v[N]$ is the quantized enveloping algebra $U_v$ generated by $E_i^{[r]}$ ($i\in I$,
$r\in\mathbb{Z}_{\ge0}$) also subject to the quantum Serre relations:
\begin{gather*}
\sum\limits_{r=0}^{1-a_{ij}} (-1)^r E_i^{[r]}E_jE_i^{[1-a_{ij}-r]}=0
\qquad
\text{for any} \ \ i\ne j\in I,
\end{gather*}
where by def\/inition $[r]_i^!E_i^{[r]}=E_i^r$.
The algebra $U_v$ is also $\mathcal{Q}$-graded by the rule $|E_i^{[r]}|=r\alpha_i$ for $i\in I$.

A powerful tool for studying certain distinguished elements of $\mathcal{A}_v[N]$ are the twisted derivations
$\theta_i,\theta_i^*:\mathcal{A}_v[N]\to\mathcal{A}_v[N]$ which satisfy $\theta_i(x_j)=\theta^*_i(x_j)=\delta_{ij}$ and
\begin{gather*}
\theta_i(xy)
=v^{(\alpha_i,|y|)}\theta_i(x)y+x\theta_i(y),
\qquad
\theta_i^*(xy) =\theta^*_i(x)y+v^{(\alpha_i,|x|)}x\theta^*_i(y).
\end{gather*}
It is well known that the maps $\varphi,\varphi^*:U_v\to\operatorname{End}(\mathcal{A}_v[N])$ given~by
$E_i\mapsto\theta_i$ or $E_i\mapsto\theta^*_i$ def\/ine respectively left and right actions of $U_v$ on
$\mathcal{A}_v[N]$.

Consider a~word $\mathbf{i}=(i_1,\ldots,i_m)\in I^m$ and def\/ine the~$m$-dimensional quantum polynomial ring
\begin{gather*}
P_\mathbf{i}=\mathbb{Z}[v^{\pm\frac{1}{2}}]\langle t_1,\ldots,t_m: t_\ell t_k=v^{(\alpha_{i_k},\alpha_{i_\ell})}t_k t_\ell
\ \text{for} \ k<\ell\rangle.
\end{gather*}
For $\mathbf{a}=(a_1,\ldots,a_m)\in\mathbb{Z}^m$ def\/ine the bar-invariant basis element $t^\mathbf{a}\in P_\mathbf{i}$
by
\begin{gather*}
t^\mathbf{a}=v^{\frac{1}{2}\sum\limits_{k<\ell}a_k a_\ell(\alpha_{i_k},\alpha_{i_\ell})}t_1^{a_1}\cdots t_m^{a_m}.
\end{gather*}

\begin{Theorem}[\cite{berenstein}] For $\mathbf{i}\in I^m$ the map $\Psi_\mathbf{i}: \mathcal{A}_v[N]\to P_\mathbf{i}$, given on
generators by $\Psi_\mathbf{i}(x_j)=\sum\limits_{k: i_k=j}t_k$ for $j\in I$, defines an algebra homomorphism.
\end{Theorem}
We will call $\Psi_\mathbf{i}$ the \emph{Feigin homomorphism} of type $\mathbf{i}$.
These homomorphisms were proposed by Boris Feigin (hence the name) as a~tool for studying the skew-f\/ield of fractions of $\mathcal{A}_v[N]$.
The Feigin homomorphisms were extensively studied in~\cite{berenstein, iohara-malikov,joseph} where many important
properties were discovered.
In particular we emphasize the following result.
\begin{Theorem}[\cite{berenstein}]
\label{th:ideal}
Suppose $\mathbf{i}\!\in\! I^m$ is a~reduced word for a~Weyl group element \mbox{$w=s_{i_1}\cdots s_{i_m}{\in} W$}.
Then the kernel $K_w:=\operatorname{ker}\Psi_\mathbf{i}$ does not depend on the choice of reduced word~$\mathbf{i}$
for~$w$.
Moreover, $K_w$ is equal to the structural ideal of the quantized coordinate ring of the closure of the unipotent cell
$N^w:=N\cap B_-wB_-$, where $B_-$ is the standard negative Borel subgroup of~$G$, i.e.~under the specialization
$v\mapsto 1$ the algebra $\mathcal{A}_v[N]/K_w$ becomes the $($integral form of$)$ the coordinate ring of $\overline{N^w}$.
\end{Theorem}

The ring $\mathcal{A}_v[\overline{N^w}]$ has a~natural choice of ``coef\/f\/icients'' which become monomials under
$\Psi_\mathbf{i}$ (see~\cite[Section~6]{berenstein-rupel} for more details), in particular they provide an Ore set by
which we may localize to obtain a~quantization of the coordinate ring of $N^w$.
In particular, this shows that the Feigin homomorphism $\Psi_\mathbf{i}$ allows one to replace the complicated algebra
$\mathcal{A}_v[N^w]$ by a~relatively simple subalgebra of quantum Laurent polynomials.
It turns out that $\Psi_\mathbf{i}$ also provides a~powerful tool for revealing quantum cluster algebra structures (see the next section for details).

\section{Quantum cluster algebras}\label{sec:qca}

A~quantum cluster algebra is a~certain type of non-commutative algebra recursively def\/ined from some initial combinatorial data.
In this section we recall the combinatorial construction of quantum cluster algebras from an initial seed.

Let~$Q$ be an acyclic quiver with vertex set~$I$ and recall the symmetrizing matrix $D=\operatorname{diag}(\mathbf{d})$,
where $\mathbf{d}=(d_i: i\in I)$.
Write $n_{ij}=n_{ji}$ for the number of arrows connecting vertices $i,j\in I$ in~$Q$ and note that all such arrows point in the same direction.
We assume further that there are only f\/initely many arrows whose source or target is the vertex $i\in I$.
From the pair $(Q,\mathbf{d})$, called a~\emph{valued quiver}, we def\/ine an adjacency matrix $B=B_Q=(b_{ij})$ by the rule
\begin{gather*}
b_{ij}=
\begin{cases}
n_{ij}d_j/\gcd(d_i,d_j) & \text{if $i\to j$ in~$Q$,}
\\
-n_{ij}d_j/\gcd(d_i,d_j) & \text{if $j\to i$ in~$Q$,}
\\
0 & \text{if $i=j$.}
\end{cases}
\end{gather*}
Notice that the matrix~$B$ is skew-symmetrizable with skew-symmetrizing matrix~$D$, i.e.~$d_ib_{ij}=-d_jb_{ji}$ for all $i,j\in I$.
To connect with the previous section we will assume that the matrix~$B$ is related to the Cartan matrix~$A$~by
$a_{ij}=-|b_{ij}|$ for $i\ne j$.

Consider an index set $J\supset I$ and a~$J\times I$ matrix $\widetilde B=(b_{ij})$ with principal $I\times I$
submatrix~$B$.
We say that a~skew-symmetric $J\times J$ matrix $\Lambda=(\lambda_{ij})$ is \emph{compatible} with~$\widetilde B$ if
\begin{gather*}
\sum\limits_{i\in J} \lambda_{ij}b_{ik}=\delta_{jk}d_k
\qquad
\text{for all} \ \   j\in J \ \text{and} \ k\in I.
\end{gather*}
Associated to~$\Lambda$ we have the quantum torus $\mathcal{T}_{\Lambda,q}$ generated over $\mathbb{Z}[v^{\pm1}]$~by
quasi-commuting \emph{cluster variables} $\mathbf{X}=(X_i: i\in J)$ subject to the relations
\begin{gather*}
X_iX_j=q^{\lambda_{ij}}X_jX_i
\qquad
\text{for} \ \ i,j\in J.
\end{gather*}
We will call the pair $\Sigma_0=(\mathbf{X},\widetilde B)$ a~\emph{quantum seed} whenever~$\Lambda$ is compatible with
$\widetilde B$, the collection~$\mathbf{X}$ is called the \emph{cluster} and~$\widetilde B$ is called the \emph{exchange matrix}.
\begin{Remark}
It is often customary to include the commutation matrix~$\Lambda$ in the data of a~quantum seed, however the quantum
cluster $\mathbf{X}$ ``knows'' its quasi-commutation so this is slightly redundant and we omit it.
\end{Remark}

In order to def\/ine the mutation of quantum seeds we need to introduce more notation.
Choose a~total order $<$ on~$J$.
Then for $\mathbf{a}\in\mathbb{Z}^J$ we may def\/ine the \emph{bar-invariant monomial} $X^\mathbf{a}\in\mathcal{T}_{\Lambda,q}$~by
\begin{gather*}
X^\mathbf{a}=v^{\sum\limits_{i<j}a_ia_j\lambda_{ji}}\vec \prod_{i\in J}X_i^{a_i},
\end{gather*}
where $\vec\prod$ denotes the product in increasing order.
Then for $k\in I$ we are now ready to def\/ine the \emph{mutation in direction~$k$},
$\mu_k\Sigma=(\mu_k\mathbf{X},\mu_k\widetilde B)$, as follows:
\begin{itemize}\itemsep=0pt
\item $\mu_k\mathbf{X}=\mathbf{X}\setminus\{X_k\}\cup\{X'_k\}$ for $X'_k$ given~by
\begin{gather*}
X'_k=X^{\mathbf{b}^k_+-\alpha_k}+X^{\mathbf{b}^k_--\alpha_k},
\end{gather*}
where $\mathbf{b}^k_+,\mathbf{b}^k_-\in\mathbb{Z}_{\ge0}^J$ are the unique vectors satisfying
$\mathbf{b}^k_+-\mathbf{b}^k_-=\mathbf{b}^k$ is the $k^{th}$ column of~$\widetilde B$;
\item $\mu_k\widetilde B=(b'_{ij})$, where
\begin{gather*}
b'_{ij}=
\begin{cases}
-b_{ij} & \text{if $i=k$ or $j=k$,}
\\
b_{ij}+[b_{ik}]_+b_{kj}+b_{ik}[-b_{kj}]_+ & \text{otherwise.}
\end{cases}
\end{gather*}
\end{itemize}
Note that the compatibility of $\tilde B$ and~$\Lambda$ ensures that $\mu_k\mathbf{X}$ is again a~quantum cluster.
Precise formulas for the commutation matrix of $\mu_k\mathbf{X}$ are known~\cite{berenstein-zelevinsky}, however we will
not need them.

Now consider a~rooted, labelled
$I$-regular tree $\mathbb{T}$ with root vertex $t_0$, where the edges emanating from each vertex are labelled 
by distinct elements of~$I$.
We assign quantum seeds $\Sigma_t=(\mathbf{X}_t,\widetilde B_t)$ to the vertices $t\in\mathbb{T}$ so that
$\Sigma_{t_0}=\Sigma_0$ and for $t\frac{ \ k }{}t'$ in $\mathbb{T}$ the seeds $\Sigma_t$ and $\Sigma_{t'}$ are related  
by the mutation in direction~$k$.
For any two (possibly distant) vertices $t,t'\in\mathbb{T}$ we say that the seeds $\Sigma_t$ and $\Sigma_{t'}$ are
\emph{mutation equivalent}.
A~fundamental result in the theory of quantum cluster algebras is the following.
\begin{Theorem}[\protect{\cite[quantum Laurent phenomenon]{berenstein-zelevinsky}}]
For any vertex $t\in\mathbb{T}$ all cluster variables of the cluster $\mathbf{X}_t$ are contained in $\mathcal{T}_{\Lambda,q}$.
\end{Theorem}
The quantum cluster algebra is then the $\mathbb{Z}[v^{\pm1}][X_j^{\pm1}:j\in J\setminus I]$-subalgebra of
$\mathcal{T}_{\Lambda,q}$ generated by all cluster variables from all seeds mutation equivalent to $\Sigma_0$:
\begin{gather*}
\mathcal{A}_q(\mathbf{X},\widetilde B)=\mathbb{Z}[v^{\pm1}][X_j^{\pm1}:j\in J\setminus I][X_{i;t}: t\in\mathbb{T},i\in
I]\subset\mathcal{T}_{\Lambda,q}.
\end{gather*}
One of the main problems in the theory of quantum cluster algebras is to explicitly describe the Laurent expansions of
all cluster variables.
In general this is still an open and active area of research.
Under the current assumption of acyclicity some progress has been made on this problem using the representation theory
of valued quivers; we present these results at the beginning of the next section.

\section{Representations of valued quivers\\ and their Hall--Ringel
algebras}\label{sec:valued quivers}

Fix a~f\/inite f\/ield $\mathbb{F}$ and an algebraic closure $\bar{\mathbb{F}}$.
For each positive integer~$d$ write $\mathbb{F}_d$ for the degree~$d$ extension of $\mathbb{F}$ in $\bar{\mathbb{F}}$.
Note that with these conventions we may intersect $\mathbb{F}_d$ and $\mathbb{F}_{d'}$ to get that
$\mathbb{F}_{\gcd(d,d')}$ is their largest common subf\/ield.
We now def\/ine representations of $(Q,\mathbf{d})$ by assigning an $\mathbb{F}_{d_i}$-vector space to each vertex $i\in
I$ and an $\mathbb{F}_{\gcd(d_i,d_j)}$-linear map to each arrow from vertex~$i$ to vertex~$j$.
The f\/inite-dimensional representations of $(Q,\mathbf{d})$ form a~f\/initary, hereditary, Abelian category
$\operatorname{rep}_\mathbb{F}(Q,\mathbf{d})$ where kernels and cokernels are taken vertex-wise, we refer the reader  
to~\cite{rupel1} for more details.

Recall that we work under the assumption that the quiver~$Q$ contains no oriented cycles.
Then the simple representations of $(Q,\mathbf{d})$ may be labeled as $S_i$ ($i\in I$) where the only non-zero vector  
space is $\mathbb{F}_{d_i}$ at vertex~$i$.
The Grothendieck group $\mathcal{K}(Q,\mathbf{d})$ of $\operatorname{rep}_\mathbb{F}(Q,\mathbf{d})$ has
a~$\mathbb{Z}$-basis given by the classes $\alpha_i=[S_i]$ ($i\in I$) of the irreducible representations $S_i$, in this
way we identify the root lattice $\mathcal{Q}$ with $\mathcal{K}(Q,\mathbf{d})$.
Write $[V]$ for the isomorphism class of a~representation $V\in\operatorname{rep}_\mathbb{F}(Q,\mathbf{d})$ and write
$|V|$ for its \emph{dimension vector}, i.e.~the class of~$V$ in $\mathcal{K}(Q,\mathbf{d})$.
Def\/ine the \emph{height} $\operatorname{ht}(V)$ of a~representation~$V$ as the length of any composition series for~$V$.

We def\/ine the \emph{Euler--Ringel form}
$\langle\cdot,\cdot\rangle:
\mathcal{K}(Q,\mathbf{d})\times\mathcal{K}(Q,\mathbf{d})\to\mathbb{Z}$ on generators~by
\begin{gather*}
\langle\alpha_i,\alpha_j\rangle=
\begin{cases}
d_i & \text{if $i=j$,}
\\
-[d_ib_{ij}]_+ & \text{if $i\ne j$.}
\end{cases}
\end{gather*}
The identif\/ication of $\mathcal{Q}$ and $\mathcal{K}(Q,\mathbf{d})$ allows to transfer the bilinear form $(\cdot,\cdot)$
to $\mathcal{K}(Q,\mathbf{d})$.
It then follows from the def\/initions that $(\cdot,\cdot)$ is the symmetrization of the Euler--Ringel form, 
i.e.~$(\alpha_i,\alpha_j)=\langle\alpha_i,\alpha_j\rangle+\langle\alpha_j,\alpha_i\rangle$ for all $i,j\in I$.
Recall that we write $\alpha_i^\vee=\alpha_i/d_i$.
For a~dimension vector $\mathbf{v}\in\mathcal{K}(Q,\mathbf{d})$ def\/ine ${}^*\mathbf{v}\in\mathcal{K}(Q,\mathbf{d})$~by
\begin{gather*}
{}^*\mathbf{v}=\sum\limits_{i\in I}\langle\alpha_i^\vee,\mathbf{v}\rangle\alpha_i.
\end{gather*}

In~\cite{rupel1} we proposed the following def\/inition and established the additive categorif\/ication of
Theorem~\ref{th:qcc} below in a~restricted set of cases.
\begin{Definition}
For any representation $V\in\operatorname{rep}_\mathbb{F}(Q,\mathbf{d})$ we associate the \emph{quantum cluster
character} $X_V$ def\/ined~by
\begin{gather*}
X_V=\sum\limits_{\mathbf{e}\in\mathcal{K}(Q,\mathbf{d})}|\mathbb{F}|^{-\frac{1}{2}\langle\mathbf{e},|V|-\mathbf{e}\rangle}
\cdot|\operatorname{Gr}_\mathbf{e}(V)|\cdot X^{-\widetilde B\mathbf{e}-{}^*|V|}\in\mathcal{T}_{\Lambda,|\mathbb{F}|},
\end{gather*}
where $\operatorname{Gr}_\mathbf{e}(V)$ denotes the set of subrepresentations of~$V$ with dimension vector $\mathbf{e}$.
\end{Definition}
A~representation~$V$ is \emph{rigid} if $\operatorname{Ext}^1(V,V)=0$.
Our main result from~\cite{rupel2} is the following.
We also refer the reader to~\cite{qin} for an analogous result in the equally valued case, i.e.~when $d_i=1$ for all
$i\in I$.
\begin{Theorem}[\protect{\cite{qin,rupel2}}]
\label{th:qcc}
Every non-initial cluster variable of $\mathcal{A}_{|\mathbb{F}|}(\mathbf{X},\widetilde B)$ is of the
form $X_V$ for some indecomposable rigid representation~$V$ of $(Q,\mathbf{d})$.
\end{Theorem}

The multiplication formulas for quantum cluster characters involved in the proof of Theorem~\ref{th:qcc} suggest
a~relationship with dual Hall--Ringel algebras.
In the work~\cite{berenstein-rupel} we studied with Arkady Berenstein the nature of this relationship, a~review of these
results are presented in the next section.

\subsection{Generalized Feigin homomorphisms}
\label{subsec:hall_gen_feigin}

Def\/ine the \emph{Hall--Ringel bialgebra} $\mathcal{H}(Q,\mathbf{d})=\bigoplus\mathbb{C}[V]$ to be the free
$\mathcal{K}(Q,\mathbf{d})$-graded $\mathbb{C}$-vector space spanned by the isomorphism classes of representations of
$(Q,\mathbf{d})$ with grading $|[V]|=|V|$ given by dimension vector.
The structure constants of the multiplication and comultiplication count certain extensions between representations.
More precisely, for representations~$U$,~$V$, and~$W$, we def\/ine the \emph{Hall number} as the cardinality of the
following set:
\begin{gather*}
\mathcal{F}_{UW}^V=\{R\subset V: R\cong W \ \text{and} \  V/R\cong U\}.
\end{gather*}
For convenience we abbreviate $v=\sqrt{|\mathbb{F}|}$.
\begin{Theorem}[\cite{ringel2}]
The map $\mu: \mathcal{H}(Q,\mathbf{d})\otimes\mathcal{H}(Q,\mathbf{d})\to\mathcal{H}(Q,\mathbf{d})$ $($written
multiplicatively$)$ given on basis elements~by
\begin{gather*}
[U][W]=\sum\limits_{[V]}v^{\langle |U|,|W|\rangle} |\mathcal{F}_{UW}^V|\cdot[V]
\end{gather*}
defines an associative multiplication on $\mathcal{H}(Q,\mathbf{d})$.
\end{Theorem}
The tensor square $\mathcal{H}(Q,\mathbf{d})\otimes\mathcal{H}(Q,\mathbf{d})$ will be considered as an algebra via the
twisted multiplication def\/ined~by
\begin{gather*}
([U]\otimes[W])([U']\otimes[W'])=v^{(|W|,|U'|)}[U][U']\otimes[W][W'].
\end{gather*}

\begin{Theorem}[\cite{green}]\label{th:hall_comult}
The map $\Delta: \mathcal{H}(Q,\mathbf{d})\to\mathcal{H}(Q,\mathbf{d})\otimes\mathcal{H}(Q,\mathbf{d})$ given on basis
elements~by
\begin{gather*}
\Delta([V])=\sum\limits_{[U],[W]}v^{\langle |U|,|W|\rangle}
\frac{|\operatorname{Aut}(U)||\operatorname{Aut}(W)|}{|\operatorname{Aut}(V)|}\cdot\big|\mathcal{F}_{UW}^V\big|\cdot[U]\otimes[W]
\end{gather*}
defines a~coassociative comultiplication on $\mathcal{H}(Q,\mathbf{d})$ which is compatible with the multiplication on
$\mathcal{H}(Q,\mathbf{d})$ and the twisted multiplication on
$\mathcal{H}(Q,\mathbf{d})\otimes\mathcal{H}(Q,\mathbf{d})$.
\end{Theorem}

The graded \emph{dual Hall--Ringel bialgebra} $\mathcal{H}^*(Q,\mathbf{d})$ has a~basis of delta functions
$\delta_{[V]}$ with grading $|\delta_{[V]}|=-|V|$.
It will be convenient to consider a~slight rescaling $[V]^*$ of the standard dual basis given~by
\begin{gather*}
[V]^*=v^{-\frac{1}{2}\langle V,V\rangle+\frac{1}{2}\sum\limits_{i=1}^nd_iv_i}\delta_{[V]}.
\end{gather*}

For a~representation $U\in\operatorname{rep}_\mathbb{F}(Q,\mathbf{d})$ def\/ine the linear maps $\theta_U,\theta_U^*:
\mathcal{H}^*(Q,\mathbf{d})\to\mathcal{H}^*(Q,\mathbf{d})$~by
\begin{gather}
\label{eq:higher derivations}
\theta_U(\delta_{[V]})([W])=\delta_{[V]}([W][U])
\qquad
\text{and}
\qquad
\theta_U^*(\delta_{[V]})([W])=\delta_{[V]}([U][W])
\end{gather}
for each $V,W\in\operatorname{rep}_\mathbb{F}(Q,\mathbf{d})$.
More precisely, we may expand $\theta_U(\delta_{[V]})$ and $\theta_U^*(\delta_{[V]})$ in the basis of delta functions as
follows:
\begin{gather*}
\theta_U(\delta_{[V]}) =\sum\limits_{[W]}v^{\langle|W|,|U|\rangle}\big|\mathcal{F}_{WU}^V\big|\cdot\delta_{[W]}, \qquad
\theta_U^*(\delta_{[V]}) =\sum\limits_{[W]}v^{\langle|U|,|W|\rangle}\big|\mathcal{F}_{UW}^V\big|\cdot\delta_{[W]}.
\end{gather*}

The following result asserts that for any $i\in I$ the linear maps $\theta_{S_i}$ and $\theta_{S_i}^*$ are twisted
derivations on $\mathcal{H}^*(Q,\mathbf{d})$.
\begin{Lemma}
For any simple representation $S_i$ of $(Q,\mathbf{d})$ and $U,W\in\operatorname{rep}_\mathbb{F}(Q,\mathbf{d})$ we have:
\begin{gather*}
\theta_{S_i}(\delta_{[U]}\delta_{[W]}) =v^{(\alpha_i,|W|)}\theta_{S_i}(\delta_{[U]})\delta_{[W]}+\delta_{[U]}\theta_{S_i}(\delta_{[W]}),
\\
\theta_{S_i}^*(\delta_{[U]}\delta_{[W]}) =\theta^*_{S_i}(\delta_{[U]})\delta_{[W]}+v^{(\alpha_i,|U|)}\delta_{[U]}\theta^*_{S_i}(\delta_{[W]}).
\end{gather*}
\end{Lemma}

\begin{proof}
These easily follow from the description~\eqref{eq:higher derivations} of $\theta_{S_i}$ and $\theta_{S_i}^*$.
Indeed, for any $V\in\operatorname{rep}_\mathbb{F}(Q,\mathbf{d})$ we have
\begin{gather*}
\theta_{S_i}(\delta_{[U]}\delta_{[W]})([V])
 =(\delta_{[U]}\delta_{[W]})([V][S_i])=(\delta_{[U]}\otimes\delta_{[W]})\big(\Delta([V][S_i])\big)
\\
\phantom{\theta_{S_i}(\delta_{[U]}\delta_{[W]})([V])}
 =(\delta_{[U]}\otimes\delta_{[W]})\big(\Delta([V])\Delta([S_i])\big)
\\
\phantom{\theta_{S_i}(\delta_{[U]}\delta_{[W]})([V])}
 =(\delta_{[U]}\otimes\delta_{[W]})\big(\Delta([V])([S_i]\otimes1)\big)+(\delta_{[U]}\otimes\delta_{[W]})\big(\Delta([V])(1\otimes[S_i])\big).
\end{gather*}
But notice that the twisted multiplication on $\mathcal{H}(Q,\mathbf{d})\otimes\mathcal{H}(Q,\mathbf{d})$ gives
\begin{gather*}
(\delta_{[U]}\otimes\delta_{[W]})\big(\Delta([V])([S_i]\otimes1)\big)
=v^{(\alpha_i,|W|)}\big(\theta_{S_1}(\delta_{[U]})\otimes\delta_{[W]}\big)\big(\Delta([V])\big)
\end{gather*}
and that
\begin{gather*}
(\delta_{[U]}\otimes\delta_{[W]})\big(\Delta([V])(1\otimes[S_i])\big)=\big(\delta_{[U]}\otimes\theta_{S_1}(\delta_{[W]})\big)\big(\Delta([V])\big).
\end{gather*}
Combining these observations completes the proof of the f\/irst identity:
\begin{gather*}
\theta_{S_i}(\delta_{[U]}\delta_{[W]})([V])=v^{(\alpha_i,|W|)}\big(\theta_{S_1}(\delta_{[U]})\otimes\delta_{[W]}\big)\big(\Delta([V])\big)
+\big(\delta_{[U]}\otimes\theta_{S_1}(\delta_{[W]})\big)\big(\Delta([V])\big)
\\
\phantom{\theta_{S_i}(\delta_{[U]}\delta_{[W]})([V])}
 =\big(v^{(\alpha_i,|W|)}\theta_{S_i}(\delta_{[U]})\delta_{[W]}+\delta_{[U]}\theta_{S_i}(\delta_{[W]})\big)([V]).
\end{gather*}
The second identity follows by a~similar calculation.
\end{proof}

\begin{Remark}
In general, one cannot hope for $\theta_U$ and $\theta_U^*$ to be derivations.
Under certain restrictions one can show that they are higher order dif\/ferential operators on
$\mathcal{H}^*(Q,\mathbf{d})$, though we will not use this here and thus leave the details for another day.
\end{Remark}

For each representation $V\in\operatorname{rep}_\mathbb{F}(Q,\mathbf{d})$ we def\/ine the \emph{quantum $\mathbf{i}$-character}
\begin{gather*}
X_{V,\mathbf{i}}=\sum\limits_{\mathbf{a}\in\mathbb{Z}_{\ge0}^m}v^{-\sum\limits_{k<\ell}
a_ka_\ell\langle\alpha_{i_\ell},\alpha_{i_k}\rangle}\cdot |\mathcal{F}_{\mathbf{i},\mathbf{a}}(V)|\cdot t^\mathbf{a}\in P_\mathbf{i},
\end{gather*}
where $\mathcal{F}_{\mathbf{i},\mathbf{a}}(V)$ denotes the set of all f\/lags of~$V$ of type $(\mathbf{i},\mathbf{a})$,
i.e.
\begin{gather*}
\mathcal{F}_{\mathbf{i},\mathbf{a}}(V)=\{0=V_m\subset V_{m-1}\subset\cdots\subset V_1\subset V_0=V: V_{k-1}/V_k\cong
S_{i_k}^{a_k},\; 1\le k\le m\}.
\end{gather*}

\begin{Theorem}[\cite{berenstein-rupel}]
\label{th:feigin_character}
For any sequence $\mathbf{i}$ the assignment $[V]^*\mapsto X_{V,\mathbf{i}}$ defines an algebra homomorphism
$\widetilde\Psi_\mathbf{i}: \mathcal{H}^*(Q,\mathbf{d})\to P_\mathbf{i}$.
\end{Theorem}

It is a~classical result of Ringel~\cite{ringel2} that the composition subalgebra of $\mathcal{H}^*(Q,\mathbf{d})$
generated over $\mathcal{A}_v=\mathbb{Z}[v^{\pm1}]$ by the classes $[S_j]^*$ of simple representations is isomorphic to
$\mathcal{A}_v[N]$.
Since $\tilde\Psi([S_j]^*)=\sum\limits_{k:i_k=j}t_k$ we see that the restriction of $\widetilde\Psi_\mathbf{i}$ to
$\mathcal{A}_v[N]$ recovers the Feigin homomor\-phism~$\Psi_\mathbf{i}$ introduced in Section~\ref{sec:feigin}.
Thus we call the maps $\widetilde\Psi_\mathbf{i}$ \emph{generalized Feigin homomorphisms}.

A certain choice of generalized Feigin homomorphism realizes the connection between acyclic quantum cluster characters
and dual Hall--Ringel algebras alluded to above.
\begin{Theorem}[\cite{berenstein-rupel}]
\label{th:char}
Let $\mathbf{i}$ be twice a~complete source adapted sequence for~$Q$.
Then there exists a~compatible pair $(\widetilde B_\mathbf{i},\Lambda_\mathbf{i})$ so that
$P_\mathbf{i}\cong\mathcal{T}_{\Lambda_\mathbf{i},q}$ and $\widetilde\Psi_\mathbf{i}([V]^*)=X_V$ for any representation
$V\in\operatorname{rep}_\mathbb{F}(Q,\mathbf{d})$.
\end{Theorem}

Combining Theorems~\ref{th:qcc} and~\ref{th:char} we see that non-initial cluster variables can be recovered as images
of rigid representations under the generalized Feigin homomorphism.

\begin{Theorem}[\cite{berenstein-rupel}]
\label{th:quantum cluster structure}
Let~$c$ be any Coxeter element of~$W$ and $\mathbf{i}_0$ a~reduced word for~$c$.
The quantized coordinate ring $\mathcal{A}_v[N^{c^2}]$ admits the structure of an acyclic quantum cluster algebra
associated to the compatible pair $(\widetilde B_\mathbf{i},\Lambda_\mathbf{i})$ for
$\mathbf{i}=(\mathbf{i}_0,\mathbf{i}_0)$.
\end{Theorem}

Together these results motivate our approach to the dual canonical basis conjecture, we will make this precise in future
sections.

\section{Quantum shuf\/f\/le algebras}\label{sec:shuffle}

The quantized coordinate ring $\mathcal{A}_v[N]$ admits an embedding into a~quantum shuf\/f\/le algebra associated to the
$I\times I$ symmetrizable Cartan matrix~$A$ with symmetrizers $\mathbf{d}$.
In this section we introduce our choice of quantum shuf\/f\/le algebra and prove some basic properties.
Better proofs of many of these results can be found in~\cite{leclerc}.
We refer the reader also to \cite{Rosso} for a more general discussion of quantum shuf\/f\/le algebras.

Let $\mathbf{W}=\bigsqcup\limits_{m\ge0} I^m$ denote the set of all words in the alphabet~$I$.
For a~word $\mathbf{i}=(i_1,\ldots,i_m)$ and a~sequence of nonnegative integers $\mathbf{c}\in\mathbb{Z}_{\ge0}^m$ we
will make use of the following simplifying notation:
\begin{gather*}
\mathbf{i}^\mathbf{c}=(i_1^{c_1},\ldots,i_m^{c_m})=(\underbrace{i_1,\ldots,
i_1}_{c_1},\ldots,\underbrace{i_m,\ldots,i_m}_{c_m}),
\end{gather*}
where $i_k$ appears $c_k$ times ($1\le k\le m$).
For $\mathbf{j}=(j_1,\ldots,j_r)\in\mathbf{W}$ we will write $|\mathbf{j}|=\alpha_{j_1}+\dots+\alpha_{j_r}$ for the
\emph{degree} of the word $\mathbf{j}$.
Denote by $\mathbf{W}^\nu$ the set of all words of degree $\nu\in\mathcal{Q}$.

Def\/ine the quantum shuf\/f\/le algebra $g^*$ as the $\mathbb{Z}[v^{\pm\frac{1}{2}}]$-algebra with basis $\mathbf{W}$ and
multiplication given~by
\begin{gather}
\label{eq:shuffle_product}
(j_1,\ldots,j_r)\circ(j_{r+1},\ldots,j_{r+s})=\sum\limits_{\sigma\in\Sigma_{r,s}}v^{\zeta(\sigma)}\mathbf{j}_\sigma,
\end{gather}
where $\mathbf{j}_\sigma=(j_{\sigma_1},\ldots,j_{\sigma_{r+s}})$ and $\zeta(\sigma)=\zeta_{r,s}(\sigma)$ measures both
inversions and non-inversions occurring in~$\sigma$:
\begin{gather*}
\zeta(\sigma)=\frac{1}{2}\sum\limits_{\substack{1\le k\le r
\\
r+1\le \ell\le r+s
\\
\sigma^{-1}_\ell<\sigma^{-1}_k}}(\alpha_{j_k},\alpha_{j_\ell})-\frac{1}{2}\sum\limits_{\substack{1\le k\le r
\\
r+1\le \ell\le r+s
\\
\sigma^{-1}_k<\sigma^{-1}_\ell}}(\alpha_{j_k},\alpha_{j_\ell}).
\end{gather*}

\begin{Remark}
Kleshchev--Ram~\cite{kleshchev-ram1} and Leclerc~\cite{leclerc} use opposite multiplications for their quantum shuf\/f\/le
algebras measuring only inversions and only non-inversions respectively.
The conventions of~\cite{kleshchev-ram1} are claimed to be more convenient for the relationship with the representation
theory of quiver Hecke algebras.

We hope that our choice of multiplication on $g^*$ will prove equally useful in this regard since this is our preferred
quantum shuf\/f\/le multiplication which will f\/igure prominently in the results presented in this note.
\end{Remark}

The following proposition can be obtained from the results of~\cite{leclerc} by rescaling the basis, however for
completeness we provide a~proof.

\begin{Proposition}
The shuffle product~\eqref{eq:shuffle_product} defines an associative algebra structure on $g^*$.
\end{Proposition}

\begin{proof}
It will be convenient to introduce the subgroup $\Sigma_{r,s,t}$ of $\Sigma_{r+s+t}$ consisting of all shuf\/f\/les of the
sequences $(1,\ldots,r)$, $(r+1,\ldots,r+s)$, and $(r+s+1,\ldots,r+s+t)$, that is
\begin{gather*}
\Sigma_{r,s,t}=\big\{\sigma\in\Sigma_{r+s+t}: \sigma^{-1}_1<\dots<\sigma^{-1}_r, \sigma^{-1}_{r+1}<\dots<\sigma^{-1}_{r+s},
\sigma^{-1}_{r+s+1}<\dots<\sigma^{-1}_{r+s+t}\big\}.
\end{gather*}
There are natural embeddings $\Sigma_{r,s}\hookrightarrow\Sigma_{r,s,t}$ written $\sigma\mapsto\widetilde\sigma$ and
$\Sigma_{s,t}\hookrightarrow\Sigma_{r,s,t}$ written $\sigma\mapsto\widetilde\sigma'$, where
\begin{gather*}
\widetilde\sigma_k=\sigma_k
\quad
\text{for}
\quad
1\le k\le r+s
\qquad
\text{and}
\qquad
\widetilde\sigma_{r+s+k}=r+s+k
\quad
\text{for}
\quad
1\le k\le t,
\\
\widetilde\sigma'_k=k
\quad
\text{for}
\quad
1\le k\le r
\qquad
\text{and}
\qquad
\widetilde\sigma'_{r+k}=r+\sigma_k
\quad
\text{for}
\quad
1\le k\le s+t.
\end{gather*}
Notice that every element of $\Sigma_{r,s,t}$ has a~unique presentation as $\widetilde\sigma\cdot\tau$ with
$\sigma\in\Sigma_{r,s}$ and $\tau\in\Sigma_{r+s,t}$ and a~unique presentation as $\widetilde\sigma'\cdot\rho$ with
$\sigma\in\Sigma_{s,t}$ and $\rho\in\Sigma_{r,s+t}$.
For $\sigma\in\Sigma_{r,s,t}$ def\/ine
\begin{gather*}
\widetilde\zeta(\sigma)=\frac{1}{2}\sum\limits_{\substack{1\le k\le r
\\
r+1\le\ell\le r+s
\\
\sigma_\ell^{-1}<\sigma_k^{-1}}}(\alpha_{j_k},\alpha_{j_\ell}) +\frac{1}{2}\sum\limits_{\substack{1\le k\le r
\\
r+s+1\le\ell\le r+s+t
\\
\sigma_\ell^{-1}<\sigma_k^{-1}}}(\alpha_{j_k},\alpha_{j_\ell}) +\frac{1}{2}\sum\limits_{\substack{r+1\le k\le r+s
\\
r+s+1\le\ell\le r+s+t
\\
\sigma_\ell^{-1}<\sigma_k^{-1}}}(\alpha_{j_k},\alpha_{j_\ell})
\\
\phantom{\widetilde\zeta(\sigma)=}{}
 -\frac{1}{2}\sum\limits_{\substack{1\le k\le r
\\
r+1\le\ell\le r+s
\\
\sigma_k^{-1}<\sigma_\ell^{-1}}}(\alpha_{j_k},\alpha_{j_\ell}) -\frac{1}{2}\sum\limits_{\substack{1\le k\le r
\\
r+s+1\le\ell\le r+s+t
\\
\sigma_k^{-1}<\sigma_\ell^{-1}}}(\alpha_{j_k},\alpha_{j_\ell}) -\frac{1}{2}\sum\limits_{\substack{r+1\le k\le r+s
\\
r+s+1\le\ell\le r+s+t
\\
\sigma_k^{-1}<\sigma_\ell^{-1}}}(\alpha_{j_k},\alpha_{j_\ell}).
\end{gather*}
Then one immediately checks that:
\begin{gather*}
\widetilde\zeta(\widetilde\sigma\cdot\tau)=\zeta(\sigma)+\zeta(\tau)
\qquad
\text{for any $\sigma\in\Sigma_{r,s}$ and $\tau\in\Sigma_{r+s,t}$,}
\\
\widetilde\zeta(\widetilde\sigma'\cdot\rho)=\zeta(\sigma)+\zeta(\rho)
\qquad
\text{for any $\sigma\in\Sigma_{s,t}$ and $\rho\in\Sigma_{r,s+t}$.}
\end{gather*}
We are now ready to prove the proposition.
Indeed, on the one hand we have
\begin{gather*}
\big((j_1,\ldots,j_r)\circ(j_{r+1},\ldots,j_{r+s})\big)\circ(j_{r+s+1},\ldots,j_{r+s+t})
\\
\qquad =\sum\limits_{\sigma\in\Sigma_{r,s}}v^{\zeta(\sigma)}\mathbf{j}_{\sigma}\circ(j_{r+s+1},\ldots,j_{r+s+t})
 =\sum\limits_{\substack{\sigma\in\Sigma_{r,s}\\ \tau\in\Sigma_{r+s,t}}}v^{\zeta(\sigma)+\zeta(\tau)}\mathbf{j}_{\widetilde\sigma\cdot\tau}
\end{gather*}
and on the other hand we have
\begin{gather*}
(j_1,\ldots,j_r)\circ\big((j_{r+1},\ldots,j_{r+s})\circ(j_{r+s+1},\ldots,j_{r+s+t})\big)
\\
\qquad =\sum\limits_{\sigma\in\Sigma_{s,t}}v^{\zeta(\sigma)}(j_1,\ldots,j_r)\circ\mathbf{j}_\sigma
 =\sum\limits_{\substack{\sigma\in\Sigma_{s,t}\\ \rho\in\Sigma_{r,s+t}}}v^{\zeta(\sigma)+\zeta(\rho)}\mathbf{j}_{\widetilde\sigma'\cdot\rho}.
\end{gather*}
But both of these expressions are equal to
$\sum\limits_{\sigma\in\Sigma_{r,s,t}}v^{\widetilde\zeta(\sigma)}\mathbf{j}_\sigma$.
\end{proof}

We will also need a~comultiplication on $g^*$.
This is def\/ined on a~word $\mathbf{h}\in\mathbf{W}$~by
\begin{gather*}
\Delta(\mathbf{h})=\sum\limits_{\mathbf{j}_1,\mathbf{j}_2\in\mathbf{W}:(\mathbf{j}_1,\mathbf{j}_2)=\mathbf{h}}
v^{\frac{1}{2}(|\mathbf{j}_1|,|\mathbf{j}_2|)}\mathbf{j}_1\otimes\mathbf{j}_2.
\end{gather*}
In the following result we consider $g^*\otimes g^*$ with the twisted multiplication
\begin{gather*}
(\mathbf{j}_{1,1}\otimes\mathbf{j}_{1,2})\circ(\mathbf{j}_{2,1}\otimes\mathbf{j}_{2,2})
=v^{(|\mathbf{j}_{1,2}|,|\mathbf{j}_{2,1}|)}(\mathbf{j}_{1,1}\circ\mathbf{j}_{2,1})\otimes(\mathbf{j}_{1,2}\circ\mathbf{j}_{2,2})
\end{gather*}
for all words $\mathbf{j}_{1,1}$, $\mathbf{j}_{1,2}$, $\mathbf{j}_{2,1}$, $\mathbf{j}_{2,2}\in\mathbf{W}$. 
\begin{Lemma}
The map $\Delta:g^*\to g^*\otimes g^*$ defines a~coassociative comultiplication on $g^*$.
Moreover, $g^*$ is a~bialgebra with counit $\epsilon:g^*\to\mathbb{C}$ given by $\varepsilon(\mathbf{j})=0$ whenever
$\mathbf{j}\ne()$ is not the empty word and $\varepsilon\big(()\big)=1$.
\end{Lemma}

\begin{proof}
We will check explicitly the coassociativity and the compatibility of multiplication and comultiplication on $g^*$, the
other axioms of a~bialgebra structure on $g^*$ are straightforward and left to the reader.

We start with the coassociativity of~$\Delta$.
For a~word $\mathbf{h}\in\mathbf{W}$ we have
\begin{gather*}
(\Delta\otimes\operatorname{Id})\Delta(\mathbf{h})
=\sum\limits_{\mathbf{j}_1,\mathbf{j}_2\in\mathbf{W}:(\mathbf{j}_1,\mathbf{j}_2)=\mathbf{h}}
v^{\frac{1}{2}(|\mathbf{j}_1|,|\mathbf{j}_2|)}\Delta(\mathbf{j}_1)
\otimes\mathbf{j}_2
\\
\phantom{(\Delta\otimes\operatorname{Id})\Delta(\mathbf{h})}
=\sum\limits_{\substack{\mathbf{j}_1,\mathbf{j}_2\in\mathbf{W}:\\(\mathbf{j}_1,\mathbf{j}_2)=\mathbf{h}}}
v^{\frac{1}{2}(|\mathbf{j}_1|,|\mathbf{j}_2|)}
\sum\limits_{\substack{\mathbf{j}_{1,1},\mathbf{j}_{1,2}\in\mathbf{W}:\\(\mathbf{j}_{1,1},\mathbf{j}_{1,2})=\mathbf{j}_1}}
v^{\frac{1}{2}(|\mathbf{j}_{1,1}|,|\mathbf{j}_{1,2}|)}(\mathbf{j}_{1,1}\otimes\mathbf{j}_{1,2})
\otimes\mathbf{j}_2
\\
\phantom{(\Delta\otimes\operatorname{Id})\Delta(\mathbf{h})}
=\sum\limits_{\mathbf{j}_{1,1},\mathbf{j}_{1,2},\mathbf{j}_2\in\mathbf{W}:(\mathbf{j}_{1,1},\mathbf{j}_{1,2},\mathbf{j}_2)=\mathbf{h}}
v^{\frac{1}{2}(|\mathbf{j}_{1,1}|+|\mathbf{j}_{1,2}|,|\mathbf{j}_2|)}
v^{\frac{1}{2}(|\mathbf{j}_{1,1}|,|\mathbf{j}_{1,2}|)}\mathbf{j}_{1,1}\otimes\mathbf{j}_{1,2}\otimes\mathbf{j}_2,
\\
(\operatorname{Id}\otimes\Delta)\Delta(\mathbf{h})
=\sum\limits_{\mathbf{j}_1,\mathbf{j}_2\in\mathbf{W}:(\mathbf{j}_1,\mathbf{j}_2)=\mathbf{h}}
v^{\frac{1}{2}(|\mathbf{j}_1|,|\mathbf{j}_2|)}\mathbf{j}_1\otimes\Delta(\mathbf{j}_2)
\\
\phantom{(\operatorname{Id}\otimes\Delta)\Delta(\mathbf{h})}
=\sum\limits_{\substack{\mathbf{j}_1,\mathbf{j}_2\in\mathbf{W}:\\(\mathbf{j}_1,\mathbf{j}_2)=\mathbf{h}}}
v^{\frac{1}{2}(|\mathbf{j}_1|,|\mathbf{j}_2|)}
\sum\limits_{\substack{\mathbf{j}_{2,1},\mathbf{j}_{2,2}\in\mathbf{W}:\\(\mathbf{j}_{2,1},\mathbf{j}_{2,2})=\mathbf{j}_2}}
v^{\frac{1}{2}(|\mathbf{j}_{2,1}|,|\mathbf{j}_{2,2}|)}\mathbf{j}_1\otimes(\mathbf{j}_{2,1}\otimes\mathbf{j}_{2,2})
\\
\phantom{(\operatorname{Id}\otimes\Delta)\Delta(\mathbf{h})}
=\sum\limits_{\mathbf{j}_1,\mathbf{j}_{2,1},\mathbf{j}_{2,2},\in\mathbf{W}:(\mathbf{j}_1,\mathbf{j}_{2,1},\mathbf{j}_{2,2})=\mathbf{h}}
v^{\frac{1}{2}(|\mathbf{j}_1|,|\mathbf{j}_{2,1}|+|\mathbf{j}_{2,2}|)}v^{\frac{1}{2}(|\mathbf{j}_{2,1}|,|\mathbf{j}_{2,2}|)}
\mathbf{j}_1\otimes\mathbf{j}_{2,1}\otimes\mathbf{j}_{2,2}.
\end{gather*}
But both of these are equal to
\begin{gather*}
\sum\limits_{\mathbf{j}_1,\mathbf{j}_2,\mathbf{j}_3,\in\mathbf{W}:(\mathbf{j}_1,\mathbf{j}_2,\mathbf{j}_3)=\mathbf{h}}
v^{\frac{1}{2}(|\mathbf{j}_1|,|\mathbf{j}_2|)}v^{\frac{1}{2}(|\mathbf{j}_1|,|\mathbf{j}_3|)}
v^{\frac{1}{2}(|\mathbf{j}_2|,|\mathbf{j}_3|)}\mathbf{j}_1\otimes\mathbf{j}_2\otimes\mathbf{j}_3,
\end{gather*}
this establishes coassociativity.

To see the compatibility of the shuf\/f\/le multiplication on $g^*$ and~$\Delta$ we consider two words
$\mathbf{j}_1,\mathbf{j}_2\in\mathbf{W}$.
On one hand we have
\begin{gather*}
\Delta(\mathbf{j}_1\circ\mathbf{j}_2)=\sum\limits_{\sigma\in\Sigma_{r,s}}v^{\zeta(\sigma)}\Delta(\mathbf{j}_\sigma)
=\sum\limits_{\sigma\in\Sigma_{r,s}}v^{\zeta(\sigma)}
\sum\limits_{\mathbf{j}_{\sigma,1},\mathbf{j}_{\sigma,2}\in\mathbf{W}:(\mathbf{j}_{\sigma,1},\mathbf{j}_{\sigma,2})=\mathbf{j}_\sigma}
v^{\frac{1}{2}(|\mathbf{j}_{\sigma,1}|,|\mathbf{j}_{\sigma,2}|)}\mathbf{j}_{\sigma,1}\otimes\mathbf{j}_{\sigma,2}.
\end{gather*}
On the other hand we have
\begin{gather*}
\Delta(\mathbf{j}_1)\circ\Delta(\mathbf{j}_2)
=\left(\sum\limits_{\substack{\mathbf{j}_{1,1},\mathbf{j}_{1,2}\in\mathbf{W}:\\(\mathbf{j}_{1,1},\mathbf{j}_{1,2})=\mathbf{j}_1}}
v^{\frac{1}{2}(|\mathbf{j}_{1,1}|,|\mathbf{j}_{1,2}|)}\mathbf{j}_{1,1}\otimes\mathbf{j}_{1,2}\right)\!
\circ\!\left(\sum\limits_{\substack{\mathbf{j}_{2,1},\mathbf{j}_{2,2}\in\mathbf{W}:\\(\mathbf{j}_{2,1},\mathbf{j}_{2,2})=\mathbf{j}_2}}
v^{\frac{1}{2}(|\mathbf{j}_{2,1}|,|\mathbf{j}_{2,2}|)}\mathbf{j}_{2,1}\otimes\mathbf{j}_{2,2}\right)
\\
\qquad
=\sum\limits_{\substack{\mathbf{j}_{1,1},\mathbf{j}_{1,2},\mathbf{j}_{2,1},\mathbf{j}_{2,2}\in\mathbf{W}:\\
(\mathbf{j}_{1,1},\mathbf{j}_{1,2})=\mathbf{j}_1,(\mathbf{j}_{2,1},\mathbf{j}_{2,2})=\mathbf{j}_2}}
v^{\frac{1}{2}(|\mathbf{j}_{1,1}|,|\mathbf{j}_{1,2}|)}v^{\frac{1}{2}(|\mathbf{j}_{2,1}|,|\mathbf{j}_{2,2}|)}
v^{(|\mathbf{j}_{1,2}|,|\mathbf{j}_{2,1}|)}(\mathbf{j}_{1,1}\circ\mathbf{j}_{2,1})\otimes(\mathbf{j}_{1,2}\circ\mathbf{j}_{2,2})
\\
\qquad
=\sum\limits_{\substack{\mathbf{j}_{1,1},\mathbf{j}_{1,2},\mathbf{j}_{2,1},\mathbf{j}_{2,2}\in\mathbf{W}:\\
(\mathbf{j}_{1,1},\mathbf{j}_{1,2})=\mathbf{j}_1,(\mathbf{j}_{2,1},\mathbf{j}_{2,2})=\mathbf{j}_2}}
\sum\limits_{\sigma_1\in\Sigma_{r_1,s_1}}
\sum\limits_{\sigma_2\in\Sigma_{r_2,s_2}}v^{\frac{1}{2}(|\mathbf{j}_{1,1}|,|\mathbf{j}_{1,2}|)}
v^{\frac{1}{2}(|\mathbf{j}_{2,1}|,|\mathbf{j}_{2,2}|)}v^{(|\mathbf{j}_{1,2}|,|\mathbf{j}_{2,1}|)}
\\
\qquad\phantom{=}  \times
v^{\zeta(\sigma_1)}v^{\zeta(\sigma_2)}\mathbf{j}_{\sigma_1,1}\otimes\mathbf{j}_{\sigma_2,2},
\end{gather*}
where we wrote $r_i$ and $s_i$ ($i=1,2$) for the lengths of $\mathbf{j}_{i,1}$ and $\mathbf{j}_{i,2}$ respectively.
There is an obvious bijection between pairs consisting of a~shuf\/f\/le $\sigma\in\Sigma_{r,s}$ together with a~partition
$(\mathbf{j}_{\sigma,1},\mathbf{j}_{\sigma,2})$ of~$\mathbf{j}_\sigma$ and quadruples consisting of partitions of~$\mathbf{j}_1$ and~$\mathbf{j}_2$ together with shuf\/f\/les $\sigma_1\in\Sigma_{r_1,s_1}$ and
$\sigma_2\in\Sigma_{r_2,s_2}$.
Indeed, reading of\/f the contributions in $\mathbf{j}_{\sigma,1}$ and $\mathbf{j}_{\sigma,2}$ coming from $\mathbf{j}_1$
and $\mathbf{j}_2$ produce the desired partitions while the relative positions of these contributions produce the
desired shuf\/f\/les in $\Sigma_{r_1,s_1}$ and $\Sigma_{r_2,s_2}$.

Notice that under this bijection we have
\begin{gather*}
\zeta(\sigma)=\zeta(\sigma_1)+\zeta(\sigma_2)+\frac{1}{2}(|\mathbf{j}_{1,2}|,|\mathbf{j}_{2,1}|)-\frac{1}{2}(|\mathbf{j}_{1,1}|,|\mathbf{j}_{2,2}|)
\end{gather*}
and
\begin{gather*}
(|\mathbf{j}_{\sigma,1}|,|\mathbf{j}_{\sigma,2}|)
=(|\mathbf{j}_{1,1}|,|\mathbf{j}_{1,2}|)+(|\mathbf{j}_{1,1}|,|\mathbf{j}_{2,2}|)
+(|\mathbf{j}_{2,1}|,|\mathbf{j}_{1,2}|)+(|\mathbf{j}_{2,1}|,|\mathbf{j}_{2,2}|).
\end{gather*}
An easy inspection now shows that $\Delta(\mathbf{j}_1\circ\mathbf{j}_2)=\Delta(\mathbf{j}_1)\circ\Delta(\mathbf{j}_2)$
as claimed.
\end{proof}

It will be useful to understand iterated shuf\/f\/le multiplication of length one words.
\begin{Lemma}
For any word $\mathbf{j}\in\mathbf{W}$ we have
\begin{gather*}
(j_1)\circ\dots\circ(j_r)=\sum\limits_{\sigma\in\Sigma_r}v^{\eta(\sigma)}\mathbf{j}_\sigma,
\end{gather*}
where $\eta(\sigma)=\frac{1}{2}\sum\limits_{\substack{k<\ell\\ \sigma_\ell^{-1}<\sigma_k^{-1}}}(\alpha_{j_k},\alpha_{j_\ell})
-\frac{1}{2}\sum\limits_{\substack{k<\ell\\ \sigma_k^{-1}<\sigma_\ell^{-1}}}(\alpha_{j_k},\alpha_{j_\ell})$.
\end{Lemma}

\begin{proof}
We work by induction on~$r$ beginning with the case $r=2$:
\begin{gather*}
(j_1)\circ(j_2)=v^{-\frac{1}{2}(\alpha_{j_1},\alpha_{j_2})}(j_1,j_2)+v^{\frac{1}{2}(\alpha_{j_1},\alpha_{j_2})}(j_2,j_1) 
\end{gather*}
which is easily seen to verify the claim.
For $r>2$, using the inductive hypothesis for~$r$ we have
\begin{gather*}
(j_1)\circ\dots\circ(j_{r+1})=\sum\limits_{\sigma\in\Sigma_r}v^{\eta(\sigma)}\mathbf{j}_\sigma\circ(j_{r+1})
=\sum\limits_{\substack{\sigma\in\Sigma_r\\ \tau\in\Sigma_{r,1}}}v^{\eta(\sigma)+\zeta(\tau)}\mathbf{j}_{\tilde\sigma\cdot\tau}
=\sum\limits_{\sigma\in\Sigma_{r+1}}v^{\eta(\sigma)}\mathbf{j}_\sigma.  \tag*{\qed}
\end{gather*}
\renewcommand{\qed}{}
\end{proof}

\begin{Lemma}
Taking $j_\ell=i$ for all~$\ell$ in the definition of $\eta(\sigma)$ we have
$\sum\limits_{\sigma\in\Sigma_r}v^{\eta(\sigma)}=[r]_i^!$.
\end{Lemma}
\begin{proof}
For $r=1$ there is nothing to show, so we begin our induction with $r=2$:
\begin{gather*}
\sum\limits_{\sigma\in\Sigma_2}v^{\eta(\sigma)}=v^{d_i}+v^{-d_i}=[2]_i.
\end{gather*}
Now assume the result holds for~$r$, then we complete the induction by a~direct calculation:
\begin{gather*}
\sum\limits_{\sigma\in\Sigma_{r+1}}v^{\eta(\sigma)}
=\sum\limits_{1\le k\le r+1}\sum\limits_{\sigma\in\Sigma_{r+1},\sigma_k^{-1}=r+1}v^{\eta(\sigma)}
 =\sum\limits_{1\le k\le r+1}v^{(r+2-2k)d_i}\sum\limits_{\sigma\in\Sigma_r}v^{\eta(\sigma)}
\\
\phantom{\sum\limits_{\sigma\in\Sigma_{r+1}}v^{\eta(\sigma)}}
 =[r+1]_i[r]_i^!.  \tag*{\qed}
\end{gather*}
\renewcommand{\qed}{}
\end{proof}

We then get the following immediate consequence.
\begin{Corollary}
\label{cor:shuffle_divided_power}
For any~$r$ we have $(i)^{\circ[r]}=(i^r)$ where $(i)^{\circ[r]}=(i)^{\circ r}/[r]_i^!$.  
\end{Corollary}

It will also be useful to understand the shuf\/f\/le of two isotypic words, the main ingredient is the following.
\begin{Lemma}
\label{lem:shuffle_binomial}
Taking $j_\ell=i$ for all~$\ell$ in the definition of $\zeta(\sigma)$, we have
\begin{gather*}
\sum\limits_{\sigma\in\Sigma_{r,s}}v^{\zeta(\sigma)}={r+s\brack r}_i
\end{gather*}
for any nonnegative integers~$r$ and~$s$.
\end{Lemma}

\begin{proof}
We proceed by induction, the base of the induction $r=s=0$ being clear.
The induction step is given by the following use of the second quantum Pascal identity from
Lemma~\ref{lem:binomial_identities}(2a):
\begin{gather*}
\sum\limits_{\sigma\in\Sigma_{r,s}}v^{\zeta(\sigma)}
 =\sum\limits_{\sigma\in\Sigma_{r,s},\sigma_{r+s}=r}\!v^{\zeta(\sigma)}+\sum\limits_{\sigma\in\Sigma_{r,s},\sigma_{r+s}=r+s}\! v^{\zeta(\sigma)}
 =v_i^s\sum\limits_{\sigma\in\Sigma_{r-1,s}}\!v^{\zeta(\sigma)}+v_i^{-r}\sum\limits_{\sigma\in\Sigma_{r,s-1}}\!v^{\zeta(\sigma)}
\\
\phantom{\sum\limits_{\sigma\in\Sigma_{r,s}}v^{\zeta(\sigma)}}
 =v_i^s{r+s-1\brack r-1}_i+v_i^{-r}{r+s-1\brack r}_i
 ={r+s\brack r}_i.  \tag*{\qed}
\end{gather*}
\renewcommand{\qed}{}
\end{proof}

This allows us to prove the following pair of results which are essential for obtaining the embedding of the quantum
group $\mathcal{A}_v[N]$ in the quantum shuf\/f\/le algebra $g^*$.
\begin{Lemma}
\label{lem:rank 2 shuffles}
For $0\le r\le n$ we have
\begin{gather*}
(i)^{\circ[r]}\circ(j)\circ(i)^{\circ[n-r]}
\\
\qquad
=\sum\limits_{s=0}^r\sum\limits_{t=0}^{n-r}
v_i^{-s(a_{ij}+n)+ta_{ij}-\frac{1}{2} na_{ij}} {s+t\brack s}_i {n-s-t\brack r-s}_i v_i^{r(a_{ij}+s+t)}(i^{s+t},j,i^{n-s-t}).
\end{gather*}
\end{Lemma}

\begin{proof}
By Corollary~\ref{cor:shuffle_divided_power} we have
$(i)^{\circ[r]}\circ(j)\circ(i)^{\circ[n-r]}=(i^r)\circ(j)\circ(i^{n-r})$.
In the product $(i^r)\circ(j)\circ(i^{n-r})$ we denote by $r-s$ the number of the f\/irst~$r$~$i$'s that shuf\/f\/le to the
right of~$j$ and denote by~$t$ the number of the last $n-r$~$i$'s that shuf\/f\/le to the left of~$j$.
The coef\/f\/icient of $(i^s,j,i^{r-s})$ in the product $(i^r)(j)$ is $v_i^{\frac{1}{2}(r-s)a_{ij}-\frac{1}{2} sa_{ij}}$.
Using Lemma~\ref{lem:shuffle_binomial}, the coef\/f\/icient of $(i^{s+t},j,i^{n-s-t})$ in the product
$(i^s,j,i^{r-s})(i^{n-r})$ is
\begin{gather*}
v_i^{(r-s)t-s(n-r-t)} {s+t\brack s}_i {n-s-t\brack r-s}_i v_i^{\frac{1}{2} ta_{ij}-\frac{1}{2} (n-r-t)a_{ij}}.
\end{gather*}
Thus, after simplifying slightly, the coef\/f\/icient of $(i^{s+t},j,i^{n-s-t})$ in the product
$(i^r)\circ(j)\circ(i^{n-r})$ becomes
\begin{gather*}
v_i^{-s(a_{ij}+n)+ta_{ij}-\frac{1}{2} na_{ij}} {s+t\brack s}_i {n-s-t\brack r-s}_i v_i^{r(a_{ij}+s+t)}
\end{gather*}
as claimed
\end{proof}

\begin{Corollary}
\label{cor:shuffle Serre}
For any $i\ne j\in I$ the following quantum Serre relation holds in $g^*$:
\begin{gather}
\label{eq:shuffle_Serre}
\sum\limits_{r=0}^{1-a_{ij}} (-1)^r (i)^{\circ[r]}\circ(j)\circ(i)^{\circ[1-a_{ij}-r]}=0.
\end{gather}
\end{Corollary}

\begin{proof}
Using Lemma~\ref{lem:rank 2 shuffles} with $n=1-a_{ij}$ we may expand the left hand side of~\eqref{eq:shuffle_Serre} as follows.
We abbreviate $\mathbf{i}_{s,t}=(i^{s+t},j,i^{1-a_{ij}-s-t})$ and in each equality we simply change the order of summation:
\begin{gather*}
 \sum\limits_{r=0}^{1-a_{ij}} (-1)^r (i)^{\circ[r]}\circ(j)\circ(i)^{\circ[1-a_{ij}-r]}
\\
\qquad = \sum\limits_{r=0}^{1-a_{ij}}\sum\limits_{s=0}^r\sum\limits_{t=0}^{1-a_{ij}-r} (-1)^r v_i^{-s+ta_{ij}-\frac{1}{2}
(1-a_{ij})a_{ij}} {s+t\brack s}_i {1-a_{ij}-s-t\brack r-s}_i v_i^{r(a_{ij}+s+t)}\mathbf{i}_{s,t}
\\
\qquad = \sum\limits_{s=0}^{1-a_{ij}}\sum\limits_{r=s}^{1-a_{ij}}\sum\limits_{t=0}^{1-a_{ij}-r} (-1)^r
v_i^{-s+ta_{ij}-\frac{1}{2} (1-a_{ij})a_{ij}} {s+t\brack s}_i {1-a_{ij}-s-t\brack r-s}_i
v_i^{r(a_{ij}+s+t)}\mathbf{i}_{s,t}
\\
 \qquad = \sum\limits_{s=0}^{1-a_{ij}}\sum\limits_{r=0}^{1-a_{ij}-s}\sum\limits_{t=0}^r (-1)^{1-a_{ij}-r}
v_i^{-s+ta_{ij}-\frac{1}{2} (1-a_{ij})a_{ij}} {s+t\brack s}_i {1-a_{ij}-s-t\brack 1-a_{ij}-r-s}_i
\\
\qquad\phantom{=}\times v_i^{(1-a_{ij}-r)(a_{ij}+s+t)}\mathbf{i}_{s,t}
\\
\qquad = \sum\limits_{s=0}^{1-a_{ij}}\sum\limits_{t=0}^{1-a_{ij}-s}\sum\limits_{r=t}^{1-a_{ij}-s} (-1)^{1-a_{ij}-r}
v_i^{-s+ta_{ij}-\frac{1}{2} (1-a_{ij})a_{ij}} {s+t\brack s}_i {1-a_{ij}-s-t\brack 1-a_{ij}-r-s}_i
\\
\qquad\phantom{=}\times
v_i^{(1-a_{ij}-r)(a_{ij}+s+t)}\mathbf{i}_{s,t}
\\
\qquad = \sum\limits_{s=0}^{1-a_{ij}}\sum\limits_{t=0}^{1-a_{ij}-s}\sum\limits_{r=0}^{1-a_{ij}-s-t} (-1)^{r+s}
v_i^{-s+ta_{ij}-\frac{1}{2} (1-a_{ij})a_{ij}} {s+t\brack s}_i {1-a_{ij}-s-t\brack r}_i\\
\qquad\phantom{=}\times
v_i^{(r+s)(a_{ij}+s+t)}\mathbf{i}_{s,t}
\\
 \qquad = \sum\limits_{s=0}^{1-a_{ij}}\sum\limits_{t=0}^{1-a_{ij}-s} (-1)^s v_i^{-s(1-a_{ij}-s-t)+ta_{ij}-\frac{1}{2}
(1-a_{ij})a_{ij}} {s+t\brack s}_i \mathbf{i}_{s,t}
\\
\qquad\phantom{=}
\times
\sum\limits_{r=0}^{1-a_{ij}-s-t} (-1)^r v_i^{r-r(1-a_{ij}-s-t)}{1-a_{ij}-s-t\brack r}_i.
\end{gather*}
By Lemma~\ref{lem:binomial_identities}, the last inner summation is zero unless $s+t=1-a_{ij}$, thus the above becomes:
\begin{gather*}
\sum\limits_{s=0}^{1-a_{ij}} (-1)^s v_i^{\frac{1}{2}(1-a_{ij})a_{ij}-sa_{ij}} {1-a_{ij} \brack s}_i \mathbf{i}_{s,1-a_{ij}-s}
\\
\qquad
=v_i^{\frac{1}{2} (1-a_{ij})a_{ij}} (i^{1-a_{ij}},j)\sum\limits_{s=0}^{1-a_{ij}}(-1)^s
v_i^{s(1-a_{ij})-s}{1-a_{ij}\brack s}_i,
\end{gather*}
which is zero by Lemma~\ref{lem:binomial_identities}.
\end{proof}

The following theorem can be proven from Corollary~\ref{cor:shuffle Serre} by the same argument as in the proof of~\cite[Theorem 4]{leclerc}.
\begin{Theorem}
There exists an injective algebra homomorphism $\kappa: \mathcal{A}_v[N]\hookrightarrow g^*$ defined on generators~by $\kappa(x_i)=(i)$.
\end{Theorem}

We now introduce certain derivations acting on $g^*$.
For $i\in I$ def\/ine maps $\theta_i,\theta_i^*: g^*\to g^*$ by linearly extending the following:
\begin{gather*}
\theta_i(\mathbf{j})=
\begin{cases}
v^{\frac{1}{2}(|\mathbf{j}|-\alpha_i,\alpha_i)}(j_1,\ldots,j_{r-1}) & \text{if $j_r=i$,}
\\
0 & \text{otherwise,}
\end{cases}
\\
\theta^*_i(\mathbf{j})=
\begin{cases}
v^{\frac{1}{2}(\alpha_i,|\mathbf{j}|-\alpha_i)}(j_2,\ldots,j_r) & \text{if $j_1=i$,}
\\
0 & \text{otherwise.}
\end{cases}
\end{gather*}

\begin{Lemma}\qquad
\begin{enumerate}[\upshape (a)]\itemsep=0pt
\item For any $i\in I$ the map $\theta_i: g^*\to g^*$ is a~twisted derivation with respect to the shuffle product, i.e.~$\theta_i$ satisfies
\begin{gather*}
\theta_i(\mathbf{j}_1\circ\mathbf{j}_2)=v^{(\alpha_i,|\mathbf{j}_2|)}\theta_i(\mathbf{j}_1)\circ\mathbf{j}_2+\mathbf{j}_1\circ\theta_i(\mathbf{j}_2)
\qquad
\text{for any}
\ \
\mathbf{j}_1,\mathbf{j}_2\in\mathbf{W}.
\end{gather*}
\item For any $i\in I$ the map $\theta_i^*: g^*\to g^*$ is a~twisted derivation with respect to the shuffle product,
i.e.~$\theta^*_i$ satisfies
\begin{gather*}
\theta^*_i(\mathbf{j}_1\circ\mathbf{j}_2)
=\theta^*_i(\mathbf{j}_1)\circ\mathbf{j}_2+v^{(|\mathbf{j}_1|,\alpha_i)}\mathbf{j}_1\circ\theta^*_i(\mathbf{j}_2)
\qquad
\text{for any}
\quad
\mathbf{j}_1,\mathbf{j}_2\in\mathbf{W}.
\end{gather*}
\end{enumerate}
\end{Lemma}

\begin{proof}
The results follow from a~direct calculation:
\begin{gather*}
\theta_i(\mathbf{j}_1\circ\mathbf{j}_2)
=\sum\limits_{\sigma\in\Sigma_{r,s}}v^{\zeta(\sigma)}\theta_i(\mathbf{j}_\sigma)
=\sum\limits_{\substack{\sigma\in\Sigma_{r,s}\\ \sigma^{-1}_r=r+s}}v^{\zeta(\sigma)}\theta_i(\mathbf{j}_\sigma)
+\sum\limits_{\substack{\sigma\in\Sigma_{r,s}\\ \sigma^{-1}_{r+s}=r+s}}v^{\zeta(\sigma)}\theta_i(\mathbf{j}_\sigma)
\\
\phantom{\theta_i(\mathbf{j}_1\circ\mathbf{j}_2)}
 =v^{(\alpha_{j_r},|\mathbf{j}_2|)}\theta_i(\mathbf{j}_1)\circ\mathbf{j}_2+\mathbf{j}_1\circ\theta_i(\mathbf{j}_2), 
\end{gather*}
and similarly:
\begin{gather*}
\theta^*_i(\mathbf{j}_1\circ\mathbf{j}_2)
=\sum\limits_{\sigma\in\Sigma_{r,s}}v^{\zeta(\sigma)}\theta^*_i(\mathbf{j}_\sigma)
=\sum\limits_{\substack{\sigma\in\Sigma_{r,s}\\ \sigma^{-1}_1=1}}v^{\zeta(\sigma)}\theta^*_i(\mathbf{j}_\sigma)
+\sum\limits_{\substack{\sigma\in\Sigma_{r,s}\\ \sigma^{-1}_{r+1}=1}}v^{\zeta(\sigma)}\theta^*_i(\mathbf{j}_\sigma)
\\
\phantom{\theta^*_i(\mathbf{j}_1\circ\mathbf{j}_2)}
 =\theta^*_i(\mathbf{j}_1)\circ\mathbf{j}_2+v^{(|\mathbf{j}_1|,\alpha_{j_{r+1}})}\mathbf{j}_1\circ\theta^*_i(\mathbf{j}_2).  \tag*{\qed}
\end{gather*}
\renewcommand{\qed}{}
\end{proof}

\subsection{Generalized Feigin homomorphisms}
\label{subsec:shuffle_gen_feigin}
Given that we have extended the Feigin homomorphisms to the dual Hall--Ringel algebra, one well-known container for the
quantum group, it is natural to try to establish extensions of the Feigin homomorphisms to the quantum shuf\/f\/le algebra.
Our proposed extension $\overline\Psi_\mathbf{i}: g^*\to P_\mathbf{i}$ agrees with the classical Feigin homomorphisms on
the generators $(i)$ ($i\in I$) of $g^*$, thus we suspect that it will play the same simplifying role as the classical
Feigin homomorphisms and our generalized Feigin homomorphisms from Section~\ref{subsec:hall_gen_feigin}.

Fix a~word $\mathbf{i}=(i_1,\ldots,i_m)\in I^m$.
We def\/ine $\overline\Psi_\mathbf{i}: g^*\to P_\mathbf{i}$ by extending linearly the following map:
\begin{gather*}
\overline\Psi_\mathbf{i}(\mathbf{j})=\sum\limits_{\substack{\mathbf{a}\in\mathbb{Z}_{\ge0}^m
\\
(i_1^{a_1},\ldots,i_m^{a_m})=\mathbf{j}}}\frac{1}{[a_1]^!_{i_1}\cdots[a_m]^!_{i_m}} t^\mathbf{a}.
\end{gather*}
Notice that $\overline\Psi_\mathbf{i}(\mathbf{j})$ is a~bar-invariant element of $P_\mathbf{i}$, hence the choice of
notation.
The following result asserts that this map def\/ines a~homomorphism of algebras.
\begin{Theorem}
The map $\overline\Psi_\mathbf{i}$ is an algebra homomorphism from $g^*$ to $P_\mathbf{i}$.
\end{Theorem}

\begin{proof}
We begin with the isotypic case $j_k=i$ ($1\le k\le r+s$) and $i_\ell=i$ ($1\le\ell\le m$).
In this case the identity
$\overline\Psi_\mathbf{i}(\mathbf{j}_1\circ\mathbf{j}_2)=\overline\Psi_\mathbf{i}(\mathbf{j}_1)\overline\Psi_\mathbf{i}(\mathbf{j}_2)$
reduces to the following:
\begin{gather*}
\sum\limits_{\sigma\in\Sigma_{r,s}}v^{\zeta(\sigma)}=\sum\limits_{\substack{\mathbf{a},\mathbf{a}'\in\mathbb{Z}_{\ge0}^m
\\
a_1+\dots+a_m=r
\\
a'_1+\dots+a'_m=s}}v_i^{\sum\limits_{k<\ell} a'_k a_\ell - \sum\limits_{k<\ell} a_k a'_\ell}{c_1\brack
a_1}_i\cdots{c_m\brack a_m}_i
\end{gather*}
for any $\mathbf{c}\in\mathbb{Z}_{\ge0}^m$ with $c_1+\dots+c_m=r+s$.
This in turn easily follows from Lemma~\ref{lem:shuffle_binomial} by repeated application of the bar-invariant subspace
identity in Lemma~\ref{lem:binomial_identities}.

Moving on to the general case we must establish the identity:
\begin{gather*}
\sum\limits_{\substack{\sigma\in\Sigma_{r,s}\\(i_1^{c_1},\ldots,i_m^{c_m})=\mathbf{j}_\sigma}}v^{\zeta(\sigma)}
=\sum\limits_{\substack{\mathbf{a},\mathbf{a}'\in\mathbb{Z}_{\ge0}^m,\mathbf{a}+\mathbf{a}'=\mathbf{c}\\
(i_1^{a_1},\ldots,i_m^{a_m})=\mathbf{j}_1\\(i_1^{a'_1},\ldots,i_m^{a'_m})=\mathbf{j}_2}}
v^{\frac{1}{2}
\sum\limits_{k<\ell} a'_k a_\ell(\alpha_{i_k},\alpha_{i_\ell})
-\frac{1}{2}\sum\limits_{\ell<k} a_\ell a'_k (\alpha_{i_k},\alpha_{i_\ell})}{c_1\brack a_1}_{i_1}\cdots{c_m\brack a_m}_{i_m}
\end{gather*}
for any $\mathbf{c}\in\mathbb{Z}_{\ge0}^m$.

For $\mathbf{a}+\mathbf{a}'=\mathbf{c}$ the term $(i_1^{c_1},\ldots,i_m^{c_m})$ in the shuf\/f\/le of
$\mathbf{j}_1=(i_1^{a_1},\ldots,i_m^{a_m})$ and $\mathbf{j}_2=(i_1^{a'_1},\ldots,i_m^{a'_m})$ is obtained by two steps:
\begin{itemize}\itemsep=0pt
\item shuf\/f\/ling the subword $i_k^{a'_k}$ of $\mathbf{j}_2$ past each subword $i_\ell^{a_\ell}$ ($\ell>k$) of
$\mathbf{j}_1$;
\item shuf\/f\/ling the subword $i_k^{a'_k}$ of $\mathbf{j}_2$ into the subword $i_k^{a_k}$ of $\mathbf{j}_1$.
\end{itemize}
The f\/irst shuf\/f\/le process produces
the exponent $\frac{1}{2}\sum\limits_{k<\ell} a'_k a_\ell (\alpha_{i_k},\alpha_{i_\ell})$ while the exponent\linebreak
$-\frac{1}{2}\sum\limits_{\ell<k} a_\ell a'_k(\alpha_{i_k},\alpha_{i_\ell})$
is a~consequence of the subword $i_k^{a'_k}$ of $\mathbf{j}_2$ not shuf\/f\/ling past any subword $i_\ell^{a_\ell}$ ($\ell<k$) of $\mathbf{j}_1$.
The second shuf\/f\/le process above produces the factor ${c_k\brack a_k}_{i_k}$ by the isotypic case.
Comparing with the right hand side of the desired equality completes the proof.
\end{proof}

\section{Quantum shuf\/f\/le characters}\label{sec:shuffle_character}

The results of Sections~\ref{sec:feigin},~\ref{subsec:shuffle_gen_feigin}, and~\ref{subsec:hall_gen_feigin} have
established the commutativity of the front two faces in the tetrahedron below:
\begin{gather}\label{dia:feigin_tetrahedron}
\begin{split}&
\includegraphics{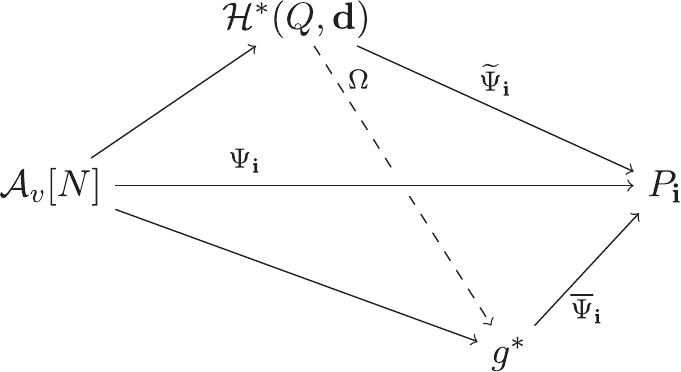}
\end{split}
\end{gather}
Our goal of this section is to complete this diagram, that is def\/ine the \emph{quantum shuffle character} homomorphism
$\Omega: \mathcal{H}^*(Q,\mathbf{d})\to g^*$ making the remaining two faces of the tetrahedron commute.

For a~word $\mathbf{j}=(j_1,\ldots,j_r)$ write $\mathcal{F}_\mathbf{j}(V)$ for the set of f\/lags in~$V$ of type~$\mathbf{j}$:
\begin{gather*}
\mathcal{F}_\mathbf{j}(V)=\{0=V_r\subset V_{r-1}\subset\dots\subset V_1\subset V_0=V: V_{k-1}/V_k\cong S_{j_k}\;\text{for $1\le k\le r$}\}. 
\end{gather*}
Then for a~basis vector $[V]^*\in\mathcal{H}^*(Q,\mathbf{d})$ we def\/ine $\Omega([V]^*)\in g^*$ as a~certain generating
function for counting f\/lags in~$V$:
\begin{gather}
\label{eq:Omega}
\Omega([V]^*)=\sum\limits_{\mathbf{j}\in\mathbf{W}} v^{-\sum\limits_{k<\ell}\langle\alpha_{j_\ell},\alpha_{j_k}\rangle}
|\mathcal{F}_\mathbf{j}(V)|\cdot\mathbf{j}.
\end{gather}
This map $\Omega:\mathcal{H}^*(Q,\mathbf{d})\to g^*$ turns out to preserve a~fair amount of structure.

\begin{Theorem}
\label{th:shuffle_character}
The map $\Omega: \mathcal{H}^*(Q,\mathbf{d})\to g^*$ given by~\eqref{eq:Omega} is a~homomorphism of twisted bialgebras.
\end{Theorem}

\begin{Remark}
In general,~$\Omega$ will not be injective.
For example, every member of the $\mathbb{P}^1$-family of indecomposable representations with dimension vector $(1,1)$
for the 2-Kronecker quiver will have the same image under~$\Omega$.
\end{Remark}

\begin{proof}
We begin by showing that~$\Omega$ is an algebra homomorphism.
Recall the rescaling
\begin{gather*}
[V]^*=v^{-\frac{1}{2}\langle |V|,|V|\rangle+\frac{1}{2}\sum\limits_{i=1}^n d_iv_i}\delta_{[V]}.
\end{gather*}
The product $[U]^*[W]^*$ can then be computed in this same basis as
\begin{gather*}
[U]^*[W]^*=v^{-\frac{1}{2}\langle |U|+|W|,|U|+|W|\rangle+\frac{1}{2}\langle |U|,|W|\rangle+\frac{1}{2}\langle
|W|,|U|\rangle+\frac{1}{2}\sum\limits_{i=1}^n d_i(u_i+w_i)}\delta_{[U]}\delta_{[W]}
\\
\phantom{[U]^*[W]^*}
 =\sum\limits_{[V]}v^{-\frac{1}{2}\langle |V|,|V|\rangle+\frac{1}{2}\langle |U|,|W|\rangle+\frac{1}{2}\langle
|W|,|U|\rangle+\frac{1}{2}\sum\limits_{i=1}^n d_iv_i}v^{\langle|U|,|W|\rangle}
\\
\phantom{[U]^*[W]^*=}{}
\times
\frac{|\operatorname{Aut}(U)||\operatorname{Aut}(W)|}{|\operatorname{Aut}(V)|}\cdot\big|\mathcal{F}^V_{UW}\big|\cdot \delta_{[V]}
\\
\phantom{[U]^*[W]^*}
 =\sum\limits_{[V]}v^{\frac{3}{2}\langle |U|,|W|\rangle+\frac{1}{2}\langle
|W|,|U|\rangle}\frac{|\operatorname{Aut}(U)||\operatorname{Aut}(W)|}{|\operatorname{Aut}(V)|}\cdot\big|\mathcal{F}^V_{UW}\big|\cdot [V]^*.
\end{gather*}
Thus the equality $\Omega([U]^*[W]^*)=\Omega([U]^*)\Omega([W]^*)$ is equivalent to the following interesting collection
of identities involving f\/lags.
\begin{Lemma}
\label{lem:flag_identity}
For every word $\mathbf{h}\in\mathbf{W}$ we have
\begin{gather}
\nonumber
\sum\limits_{[V]} v^{\frac{3}{2}\langle |U|,|W|\rangle+\frac{1}{2}\langle|W|,|U|\rangle-\sum\limits_{k<\ell}\langle\alpha_{h_\ell},\alpha_{h_k}\rangle}
\frac{|\operatorname{Aut}(U)||\operatorname{Aut}(W)|}{|\operatorname{Aut}(V)|}\big|\mathcal{F}^V_{UW}\big||\mathcal{F}_\mathbf{h}(V)|
\\
\qquad
 =\sum\limits_{\substack{\mathbf{j}_1,\mathbf{j}_2\in\mathbf{W}\\
\sigma\in\Sigma_{\operatorname{ht}(U),\operatorname{ht}(W)}: \mathbf{j}_\sigma=\mathbf{h}}}v^{-\sum\limits_{1\le
k<\ell\le t}\langle\alpha_{j_\ell},\alpha_{j_k}\rangle-\sum\limits_{t+1\le k<\ell\le
r+s}\langle\alpha_{j_\ell},\alpha_{j_k}\rangle+\zeta(\sigma)}|\mathcal{F}_{\mathbf{j}_1}(U)||\mathcal{F}_{\mathbf{j}_2}(W)|,
\label{eq:flag_identity}
\end{gather}
where we abbreviated $t=\operatorname{ht}(U)$.
\end{Lemma}

\begin{proof}
This essentially follows from the compatibility of multiplication and comultiplication in the Hall--Ringel algebra
$\mathcal{H}(Q,\mathbf{d})$.
Using the def\/initions it is easy to compute the following composition of comultiplication and iterated multiplication:
\begin{gather*}
\Delta([S_{h_1}]\cdots[S_{h_{r+s}}])=\sum\limits_{[V]}v^{\sum\limits_{k<\ell}\langle
\alpha_{h_k},\alpha_{h_\ell}\rangle}|\mathcal{F}_\mathbf{h}(V)|\Delta([V])
\\
\qquad
 =\sum\limits_{[U],[V],[W]}v^{\langle|U|,|W|\rangle+\sum\limits_{k<\ell}\langle\alpha_{h_k},\alpha_{h_\ell}\rangle}
\frac{|\operatorname{Aut}(U)||\operatorname{Aut}(W)|}{|\operatorname{Aut}(V)|}\big|\mathcal{F}^V_{UW}\big||\mathcal{F}_\mathbf{h}(V)|\cdot[U]\otimes[V].
\end{gather*}
Theorem~\ref{th:hall_comult} asserts that this is the same as the following iterated multiplication of
comultiplications:
\begin{gather*}
\Delta([S_{h_1}])\cdots\Delta([S_{h_{r+s}}])=([S_{h_1}]\otimes1+1\otimes[S_{h_1}])\cdots([S_{h_{r+s}}]\otimes1+1\otimes[S_{h_{r+s}}])
\\
\qquad
 =\sum\limits_{\substack{[U],[W]:\\ \operatorname{ht}(U)+\operatorname{ht}(W)=r+s\\
\sigma\in\Sigma_{\operatorname{ht}(U),\operatorname{ht}(W)}\\
\mathbf{j}_1,\mathbf{j}_2\in\mathbf{W}: \mathbf{j}_\sigma=\mathbf{h}}}v^E
|\mathcal{F}_{\mathbf{j}_1}(U)||\mathcal{F}_{\mathbf{j}_2}(W)|\cdot [U]\otimes[W],
\end{gather*}
with exponent $E=\sum\limits_{1\le k<\ell\le t}\langle\alpha_{j_k},\alpha_{j_\ell}\rangle
+\sum\limits_{t+1\le k<\ell\le r+s}\langle\alpha_{j_k},\alpha_{j_\ell}\rangle
+\sum\limits_{\substack{1\le k\le t\\t+1\le\ell\le r+s\\\sigma_k^{-1}>\sigma_\ell^{-1}}}(\alpha_{j_k},\alpha_{j_\ell})$,
where we abbreviate $t=\operatorname{ht}(U)$.

A few remarks seem necessary to clarify the second equality above.
To obtain a~f\/ixed $[U]$ and $[W]$ we must take from the product
\begin{gather*}
([S_{h_1}]\otimes1+1\otimes[S_{h_1}])\cdots([S_{h_{r+s}}]\otimes1+1\otimes[S_{h_{r+s}}])
\end{gather*}
exactly $\operatorname{ht}(U)$ factors, recorded by $\mathbf{j}_1$, of the form $[S_{h_k}]\otimes1$ and
$\operatorname{ht}(W)$ factors, recorded by $\mathbf{j}_2$, of the form $1\otimes[S_{h_k}]$.
With this notation we are looking at the products $[S_{j_1}]\cdots[S_{j_t}]\otimes[S_{j_{t+1}}]\cdots[S_{j_{r+s}}]$, where
again we have used the abbreviation $\operatorname{ht}(U)=t$.
The f\/irst and second sums in~$E$ as well as the f\/lag coef\/f\/icients should now be transparent.
In order to record the contribution of the twisted multiplication on
$\mathcal{H}(Q,\mathbf{d})\otimes\mathcal{H}(Q,\mathbf{d})$ to the coef\/f\/icient of $[U]\otimes[W]$ we need a~shuf\/f\/le
$\sigma\in\Sigma_{\operatorname{ht}(U),\operatorname{ht}(W)}$ and $j_k=h_{\sigma_k^{-1}}$.
Indeed, the condition $\sigma_k^{-1}>\sigma_\ell^{-1}$ of the last sum in~$E$ records a~product of the form
$(1\otimes[S_{j_k}])([S_{j_\ell}]\otimes1)$ and thus produces the stated contribution.

As previously stated, Theorem~\ref{th:hall_comult} allows us, for a~f\/ixed $[U]\otimes[W]$, to extract the equality below
from which the claim will follow:
\begin{gather}
\sum\limits_{[V]}v^{\langle|U|,|W|\rangle+\sum\limits_{k<\ell}\langle\alpha_{h_k},\alpha_{h_\ell}\rangle}
\frac{|\operatorname{Aut}(U)||\operatorname{Aut}(W)|}{|\operatorname{Aut}(V)|}|\mathcal{F}^V_{UW}||\mathcal{F}_\mathbf{h}(V)|
\nonumber
\\
\qquad
=\sum\limits_{\substack{\sigma\in\Sigma_{\operatorname{ht}(U),\operatorname{ht}(W)}\\
\mathbf{j}_1,\mathbf{j}_2\in\mathbf{W}: \mathbf{j}_\sigma=\mathbf{h}}}
v^E|\mathcal{F}_{\mathbf{j}_1}(U)||\mathcal{F}_{\mathbf{j}_2}(W)|.
\label{eq:green_equality}
\end{gather}
To complete the proof we perform a~few simple manipulations with the sums appearing in the exponents
of~\eqref{eq:flag_identity} and~\eqref{eq:green_equality}.
To start we match the left hand sides of~\eqref{eq:flag_identity} and~\eqref{eq:green_equality} by multiplying both
sides of~\eqref{eq:flag_identity} by $v^{\sum\limits_{k<\ell}\langle
\alpha_{h_\ell},\alpha_{h_k}\rangle-\frac{1}{2}\langle |U|,|W|\rangle-\frac{1}{2}\langle |W|,|U|\rangle}$ and both sides
of~\eqref{eq:green_equality} by $v^{-\sum\limits_{k<\ell}\langle \alpha_{h_k},\alpha_{h_\ell}\rangle}$.
Thus it remains to establish the following identity of exponents.
\begin{Lemma}
\label{lem:exp_equality}
For any $\sigma\in\Sigma_{\operatorname{ht}(U),\operatorname{ht}(W)}$ and words $\mathbf{j}_1$, $\mathbf{j}_2$ satisfying
$\mathcal{F}_{\mathbf{j}_1}(U)\ne\varnothing$, $\mathcal{F}_{\mathbf{j}_2}(W)\ne\varnothing$, and
$\mathbf{j}_\sigma=\mathbf{h}$ we have:
\begin{gather}
\nonumber
\sum\limits_{k<\ell}\langle \alpha_{h_\ell},\alpha_{h_k}\rangle
-\!\sum\limits_{1\le k<\ell\le t}\!\!\langle\alpha_{j_\ell},\alpha_{j_k}\rangle
-\!\sum\limits_{t+1\le k<\ell\le r+s}\!\!\langle\alpha_{j_\ell},\alpha_{j_k}\rangle
+\zeta(\sigma)-\frac{1}{2}\langle |U|,|W|\rangle-\frac{1}{2}\langle|W|,|U|\rangle
\nonumber
\\
\qquad
=-\sum\limits_{k<\ell}\langle \alpha_{h_k},\alpha_{h_\ell}\rangle+\sum\limits_{1\le k<\ell\le t}\!\!\langle\alpha_{j_k},\alpha_{j_\ell}\rangle
+\sum\limits_{t+1\le k<\ell\le r+s}\!\!\!\!\langle\alpha_{j_k},\alpha_{j_\ell}\rangle
+\!\!\!\sum\limits_{\substack{1\le k\le t\\t+1\le\ell\le r+s\\ \sigma_k^{-1}>\sigma_\ell^{-1}}}\!\!\!\!(\alpha_{j_k},\alpha_{j_\ell}). 
\label{eq:exp_equality}
\end{gather}
\end{Lemma}

\begin{proof}
We f\/irst eliminate the dependence on $\mathbf{h}$ using the following identities which are immediate consequences of the
equality $j_k=h_{\sigma_k^{-1}}$:
\begin{gather*}
\sum\limits_{k<\ell}\langle \alpha_{h_\ell},\alpha_{h_k}\rangle
-\sum\limits_{1\le k<\ell\le t}\langle\alpha_{j_\ell},\alpha_{j_k}\rangle
-\sum\limits_{t+1\le k<\ell\le r+s}\langle\alpha_{j_\ell},\alpha_{j_k}\rangle
\\
\qquad
=\sum\limits_{\substack{1\le k\le t\\t+1\le\ell\le r+s\\ \sigma_k^{-1}<\sigma_\ell^{-1}}}
\langle\alpha_{j_\ell},\alpha_{j_k}\rangle
+\sum\limits_{\substack{1\le k\le t\\t+1\le\ell\le r+s\\ \sigma_\ell^{-1}<\sigma_k^{-1}}}
\langle\alpha_{j_k},\alpha_{j_\ell}\rangle,
\\
-\sum\limits_{k<\ell}\langle \alpha_{h_k},\alpha_{h_\ell}\rangle
+\sum\limits_{1\le k<\ell\le t}\langle\alpha_{j_k},\alpha_{j_\ell}\rangle
+\sum\limits_{t+1\le k<\ell\le r+s}\langle\alpha_{j_k},\alpha_{j_\ell}\rangle
\\
\qquad
=-\sum\limits_{\substack{1\le k\le t\\t+1\le\ell\le r+s\\ \sigma_k^{-1}<\sigma_\ell^{-1}}}
\langle\alpha_{j_k},\alpha_{j_\ell}\rangle
-\sum\limits_{\substack{1\le k\le t\\t+1\le\ell\le r+s\\ \sigma_\ell^{-1}<\sigma_k^{-1}}}
\langle\alpha_{j_\ell},\alpha_{j_k}\rangle,
\end{gather*}
this reduces~\eqref{eq:exp_equality} to showing
\begin{gather}
\nonumber
\sum\limits_{\substack{1\le k\le t\\t+1\le\ell\le r+s\\ \sigma_k^{-1}<\sigma_\ell^{-1}}}\langle\alpha_{j_\ell},\alpha_{j_k}\rangle
+\sum\limits_{\substack{1\le k\le t\\t+1\le\ell\le r+s\\ \sigma_\ell^{-1}<\sigma_k^{-1}}}\langle\alpha_{j_k},\alpha_{j_\ell}\rangle
+\zeta(\sigma)-\frac{1}{2}\langle|U|,|W|\rangle-\frac{1}{2}\langle |W|,|U|\rangle
\\
\qquad
=-\sum\limits_{\substack{1\le k\le t\\t+1\le\ell\le r+s\\ \sigma_k^{-1}<\sigma_\ell^{-1}}}\langle\alpha_{j_k},\alpha_{j_\ell}\rangle
-\sum\limits_{\substack{1\le k\le t\\t+1\le\ell\le r+s\\ \sigma_\ell^{-1}<\sigma_k^{-1}}}\langle\alpha_{j_\ell},\alpha_{j_k}\rangle
+\sum\limits_{\substack{1\le k\le t\\t+1\le\ell\le r+s\\ \sigma_k^{-1}>\sigma_\ell^{-1}}}(\alpha_{j_k},\alpha_{j_\ell}).
\label{eq:exp_equality2}
\end{gather}
Using the dichotomy $\sigma_k^{-1}<\sigma_\ell^{-1}$ or $\sigma_\ell^{-1}<\sigma_k^{-1}$ for $1\le k\le t$ and
$t+1\le\ell\le r+s$ the last sum on the right hand side of~\eqref{eq:exp_equality2} becomes:
\begin{gather}
\label{eq:zeta}
\sum\limits_{\substack{1\le k\le t\\t+1\le\ell\le r+s\\ \sigma_k^{-1}>\sigma_\ell^{-1}}}(\alpha_{j_k},\alpha_{j_\ell})
=\zeta(\sigma)+\frac{1}{2}\sum\limits_{\substack{1\le k\le t\\t+1\le\ell\le r+s}}(\alpha_{j_k},\alpha_{j_\ell}).
\end{gather}
Rewriting
$(\alpha_{j_k},\alpha_{j_\ell})=\langle\alpha_{j_k},\alpha_{j_\ell}\rangle+\langle\alpha_{j_\ell},\alpha_{j_k}\rangle$
and applying~\eqref{eq:zeta} in the right hand side of~\eqref{eq:exp_equality2} gives:
\begin{gather*}
-\sum\limits_{\substack{1\le k\le t\\t+1\le\ell\le r+s\\ \sigma_k^{-1}<\sigma_\ell^{-1}}}\langle\alpha_{j_k},\alpha_{j_\ell}\rangle
-\sum\limits_{\substack{1\le k\le t\\t+1\le\ell\le r+s\\ \sigma_\ell^{-1}<\sigma_k^{-1}}}\langle\alpha_{j_\ell},\alpha_{j_k}\rangle
+\zeta(\sigma)+\frac{1}{2}\sum\limits_{\substack{1\le k\le t\\t+1\le\ell\le r+s}}
(\langle\alpha_{j_k},\alpha_{j_\ell}\rangle+\langle\alpha_{j_\ell},\alpha_{j_k}\rangle)
\\
\qquad
=-\frac{1}{2}\sum\limits_{\substack{1\le k\le t\\t+1\le\ell\le r+s}}
(\langle\alpha_{j_k},\alpha_{j_\ell}\rangle+\langle\alpha_{j_\ell},\alpha_{j_k}\rangle)+\zeta(\sigma)
\\
\qquad
\phantom{=}
+\sum\limits_{\substack{1\le k\le t\\t+1\le\ell\le r+s\\ \sigma_k^{-1}<\sigma_\ell^{-1}}}\langle\alpha_{j_\ell},\alpha_{j_k}\rangle
+\sum\limits_{\substack{1\le k\le t\\t+1\le\ell\le r+s\\ \sigma_\ell^{-1}<\sigma_k^{-1}}}\langle\alpha_{j_k},\alpha_{j_\ell}\rangle
\\
\qquad
=-\frac{1}{2}\langle |U|,|W|\rangle-\frac{1}{2}\langle |W|,|U|\rangle+\zeta(\sigma)
+\sum\limits_{\substack{1\le k\le t\\t+1\le\ell\le r+s\\ \sigma_k^{-1}<\sigma_\ell^{-1}}}\langle\alpha_{j_\ell},\alpha_{j_k}\rangle
+\sum\limits_{\substack{1\le k\le t\\t+1\le\ell\le r+s\\ \sigma_\ell^{-1}<\sigma_k^{-1}}}\langle\alpha_{j_k},\alpha_{j_\ell}\rangle
\end{gather*}
as desired, which completes the proof of Lemma~\ref{lem:exp_equality}.
\end{proof}

This completes the proof of Lemma~\ref{lem:flag_identity}.
\end{proof}

With Lemma~\ref{lem:flag_identity} proven, we may conclude that~$\Omega$ def\/ines a~homomorphism of algebras.
It remains to show that $\Delta\Omega([V]^*)=(\Omega\otimes\Omega)(\Delta([V]^*))$.
The coproduct $\Delta([V]^*)$ can be expanded in the same rescaled basis as
\begin{gather*}
\Delta([V]^*)=v^{-\frac{1}{2}\langle |V|,|V|\rangle+\frac{1}{2}\sum\limits_{i=1}^n d_iv_i}\Delta(\delta_{[V]})
\\
\phantom{\Delta([V]^*)}
 =\sum\limits_{[U],[W]}v^{-\frac{1}{2}\langle |U|+|W|,|U|+|W|\rangle+\frac{1}{2}\sum\limits_{i=1}^n
d_i(u_i+w_i)}v^{\langle|U|,|W|\rangle}\big|\mathcal{F}^V_{UW}\big|\cdot \delta_{[U]}\otimes\delta_{[W]}
\\
\phantom{\Delta([V]^*)}
 =\sum\limits_{[V]}v^{\frac{1}{2}\langle |U|,|W|\rangle-\frac{1}{2}\langle |W|,|U|\rangle} \big|\mathcal{F}^V_{UW}\big|\cdot
[U]^*\otimes[W]^*.
\end{gather*}
Then the identity $\Delta\Omega([V]^*)=(\Omega\otimes\Omega)(\Delta([V]^*))$ is equivalent to the following collection
of identities involving f\/lags.
\begin{Lemma}
For all words $\mathbf{j}_1,\mathbf{j}_2\in\mathbf{W}$ we have
\begin{gather*}
 v^{-\sum\limits_{k<\ell}\langle
\alpha_{h_\ell},\alpha_{h_k}\rangle}v^{\frac{1}{2}(|\mathbf{j}_1|,|\mathbf{j}_2|)}|\mathcal{F}_\mathbf{h}(V)|
 =\sum\limits_{[U],[W]}v^{-\sum\limits_{1\le k<\ell\le r}\langle\alpha_{j_\ell},\alpha_{j_k}\rangle}
v^{-\sum\limits_{r+1\le k<\ell\le r+s}\langle\alpha_{j_\ell},\alpha_{j_k}\rangle}
\\
\qquad
\phantom{=}{}
\times
v^{\frac{1}{2}\langle|U|,|W|\rangle-\frac{1}{2}\langle|W|,|U|\rangle}
\big|\mathcal{F}_{U,W}^V\big|\cdot|\mathcal{F}_{\mathbf{j}_1}(U)|\cdot|\mathcal{F}_{\mathbf{j}_2}(W)|,
\end{gather*}
where we used the word $\mathbf{h}=(\mathbf{j}_1,\mathbf{j}_2)$.
\end{Lemma}

\begin{proof}
This essentially follows from the associativity of multiplication in the Hall--Ringel algebra
$\mathcal{H}(Q,\mathbf{d})$.
Indeed, expanding the products
\begin{gather*}
[S_{h_1}]\cdots[S_{h_{r+s}}]=([S_{j_1}]\cdots[S_{j_r}])\cdot([S_{j_{r+1}}]\cdots[S_{j_{r+s}}])
\end{gather*}
gives
$|\mathcal{F}_\mathbf{h}(V)|=|\mathcal{F}_{U,W}^V|\cdot|\mathcal{F}_{\mathbf{j}_1}(U)|\cdot|\mathcal{F}_{\mathbf{j}_2}(W)|$.
To see the equality of the exponents on either side we note that
\begin{gather*}
\frac{1}{2}(\mathbf{j}_1,\mathbf{j}_2)
=\frac{1}{2}\langle|\mathbf{j}_1|,|\mathbf{j}_2|\rangle+\frac{1}{2}\langle|\mathbf{j}_2|,|\mathbf{j}_1|\rangle
=\frac{1}{2}\langle|U|,|W|\rangle+\frac{1}{2}\langle|W|,|U|\rangle.
\end{gather*}
Combining this with the equality
\begin{gather*}
-\sum\limits_{k<\ell}\langle \alpha_{h_\ell},\alpha_{h_k}\rangle=-\sum\limits_{1\le k<\ell\le r}\langle
\alpha_{j_\ell},\alpha_{j_k}\rangle-\sum\limits_{r+1\le k<\ell\le r+s}\langle
\alpha_{j_\ell},\alpha_{j_k}\rangle-\langle|\mathbf{j}_2|,|\mathbf{j}_1|\rangle
\end{gather*}
completes the proof.
\end{proof}

To f\/inish the proof of Theorem~\ref{th:shuffle_character} we note that~$\Omega$ respects the grading and thus, since the
multiplication on both $\mathcal{H}^*(Q,\mathbf{d})\otimes\mathcal{H}^*(Q,\mathbf{d})$ and $g^*\otimes g^*$ are twisted
by the grading, we may conclude that~$\Omega$ is a~homomorphism of twisted bialgebras.
\end{proof}

We now show that this def\/inition of~$\Omega$ does in fact make the tetrahedron commute.
The left hand triangle is easily seen to be commutative by looking at generators of $\mathcal{A}_v[N]$.
Thus we focus on the right hand triangle.

\begin{Theorem}
\label{th:psi factorization}
There is an equality of homomorphisms $\widetilde\Psi_\mathbf{i}=\overline\Psi_\mathbf{i}\circ\Omega$.
\end{Theorem}

\begin{proof}
It suf\/f\/ices to check the equality on the basis $[V]^*$ of $\mathcal{H}^*(Q,\mathbf{d})$.
Indeed, the value of $\widetilde\Psi_\mathbf{i}$ on $[V]^*$ is provided by Theorem~\ref{th:feigin_character}:
\begin{gather*}
\widetilde\Psi_\mathbf{i}([V]^*)=\sum\limits_{\mathbf{a}\in\mathbb{Z}_{\ge0}^m}v^{-\sum\limits_{k<\ell}
a_ka_\ell\langle\alpha_{i_\ell},\alpha_{i_k}\rangle}|\mathcal{F}_{\mathbf{i},\mathbf{a}}(V)|t^\mathbf{a}.
\end{gather*}
Observe that the passage from $\mathcal{F}_{\mathbf{i},\mathbf{a}}(V)$ to $\mathcal{F}_{\mathbf{i}^\mathbf{a}}(V)$
amounts to choosing a~basis in each quo\-tient~$S_{i_k}^{a_k}$ for which there are $(a_k)_{i_k}^!$ possibilities.
Thus we have the identity:
\begin{gather}
\nonumber
\widetilde\Psi_\mathbf{i}([V]^*)=\sum\limits_{\mathbf{a}\in\mathbb{Z}_{\ge0}^m}v^{-\sum\limits_{k<\ell}
a_ka_\ell\langle\alpha_{i_\ell},\alpha_{i_k}\rangle}|\mathcal{F}_{\mathbf{i}^\mathbf{a}}(V)|\frac{1}{(a_1)_{i_1}^!\cdots(a_m)_{i_m}^!}t^\mathbf{a}
\\
\label{eq:feigin_equality_LHS}
\phantom{\widetilde\Psi_\mathbf{i}([V]^*)}
 =\sum\limits_{\mathbf{a}\in\mathbb{Z}_{\ge0}^m}v^{-\sum\limits_{k<\ell}
a_ka_\ell\langle\alpha_{i_\ell},\alpha_{i_k}\rangle}|\mathcal{F}_{\mathbf{i}^\mathbf{a}}(V)|\frac{v^{-\sum\limits_{k=1}^m
d_{i_k}{a_k\choose 2}}}{[a_1]_{i_1}^!\cdots[a_m]_{i_m}^!}t^\mathbf{a},
\end{gather}
where we passed to the bar invariant~$q$-numbers in~\eqref{eq:feigin_equality_LHS}.

Now working from the def\/initions of~$\Omega$ and $\overline\Psi_\mathbf{i}$ we get:
\begin{gather}
\label{eq:feigin_equality_RHS}
\overline\Psi_\mathbf{i}\circ\Omega([V]^*)
=\sum\limits_{\mathbf{a}\in\mathbb{Z}_{\ge0}^m}v^{-\sum\limits_{k<\ell}\langle\alpha_{j_\ell},\alpha_{j_k}\rangle}
|\mathcal{F}_{\mathbf{i}^\mathbf{a}}(V)|\frac{1}{[a_1]_{i_1}^!\cdots[a_m]_{i_m}^!}t^\mathbf{a},
\end{gather}
where we have used $\mathbf{j}=\mathbf{i}^\mathbf{a}=(\underbrace{i_1,\ldots, i_1}_{a_1},\ldots,\underbrace{i_m,\ldots,i_m}_{a_m})$
in the exponent of~$v$ in~\eqref{eq:feigin_equality_RHS}.
But notice that
\begin{gather*}
\sum\limits_{k<\ell}\langle\alpha_{j_\ell},\alpha_{j_k}\rangle=d_{i_k}{a_k\choose 2}+\sum\limits_{k<\ell}
a_ka_\ell\langle\alpha_{i_\ell},\alpha_{i_k}\rangle,
\end{gather*}
where $d_{i_k}{a_k\choose 2}$ comes from taking $j_k$ and $j_\ell$ from within a~single grouping $\underbrace{i_k,\ldots, i_k}_{a_k}$
while\linebreak %
$\sum\limits_{k<\ell} a_ka_\ell\langle\alpha_{i_\ell},\alpha_{i_k}\rangle$ is obtained by taking $j_k$ and $j_\ell$ from dif\/ferent groupings.
Thus we see that $\widetilde \Psi_\mathbf{i}([V]^*)=\overline\Psi_\mathbf{i}\circ\Omega([V]^*)$ for any representation
$V\in\operatorname{rep}_\mathbb{F}(Q,\mathbf{d})$.
\end{proof}

In fact, there is an extra bit of structure preserved by the quantum shuf\/f\/le character~$\Omega$.
\begin{Proposition}
\qquad
\begin{enumerate}[\upshape (a)]\itemsep=0pt
\item The quantum shuffle character~$\Omega$ intertwines the action of $\theta_{S_i}$ on $\mathcal{H}^*(Q,\mathbf{d})$
and the action of $\theta_i$ on $g^*$, i.e.~$\theta_i\circ\Omega=\Omega\circ\theta_{S_i}$.
\item The quantum shuffle character~$\Omega$ intertwines the action of $\theta^*_{S_i}$ on $\mathcal{H}^*(Q,\mathbf{d})$
and the action of $\theta^*_i$ on $g^*$, i.e.~$\theta^*_i\circ\Omega=\Omega\circ\theta^*_{S_i}$.
\end{enumerate}
\end{Proposition}

\begin{proof}
To prove (a) we simply compute the actions on a~basis vector $[V]^*$:
\begin{gather}
\nonumber
\Omega\circ\theta_{S_i}([V]^*)=v^{-\frac{1}{2}\langle V,V\rangle+\frac{1}{2}\sum\limits_{i=1}^nd_iv_i}\Omega\circ\theta_{S_i}(\delta_{[V]})
\\
\phantom{\Omega\circ\theta_{S_i}([V]^*)}
\nonumber
=v^{-\frac{1}{2}\langle V,V\rangle+\frac{1}{2}\sum\limits_{i=1}^nd_iv_i}\sum\limits_{[W]}v^{\langle|W|,\alpha_i\rangle}
\big|\mathcal{F}_{WS_i}^V\big|\cdot\Omega(\delta_{[W]})
\\
\phantom{\Omega\circ\theta_{S_i}([V]^*)}
\nonumber
 =\sum\limits_{[W]}v^{\frac{1}{2}\langle
W,\alpha_i\rangle-\frac{1}{2}\langle\alpha_i,W\rangle}\big|\mathcal{F}_{WS_i}^V\big|\cdot\Omega([W]^*)
\\
\phantom{\Omega\circ\theta_{S_i}([V]^*)}
\nonumber
 =\sum\limits_{[W]}\sum\limits_{\mathbf{j}\in\mathbf{W}} v^{\frac{1}{2}\langle
W,\alpha_i\rangle-\frac{1}{2}\langle\alpha_i,W\rangle}v^{-\sum\limits_{k<\ell}\langle\alpha_{j_\ell},\alpha_{j_k}\rangle}
\big|\mathcal{F}_{WS_i}^V\big|\cdot|\mathcal{F}_\mathbf{j}(W)|\cdot\mathbf{j}
\\
\phantom{\Omega\circ\theta_{S_i}([V]^*)}
=\sum\limits_{\mathbf{j}\in\mathbf{W}} v^{\frac{1}{2}(|\mathbf{j}|,\alpha_i)-\langle
\alpha_i,|\mathbf{j}|\rangle}v^{-\sum\limits_{k<\ell}\langle\alpha_{j_\ell},\alpha_{j_k}\rangle}
|\mathcal{F}_{(\mathbf{j},i)}(V)|\cdot\mathbf{j},
\label{eq:derivation_RHS}
\end{gather}
where we have used in the last equality that
$\sum\limits_{[W]}|\mathcal{F}_{WS_i}^V|\cdot|\mathcal{F}_\mathbf{j}(W)|=|\mathcal{F}_{(\mathbf{j},i)}(V)|$ and
$|W|=|\mathbf{j}|$ whenever $\mathcal{F}_\mathbf{j}(W)\ne\varnothing$.
On the other hand we have:
\begin{gather}
\theta_i\circ\Omega([V]^*)
\label{eq:derivation_LHS}
 =\sum\limits_{\mathbf{j}\in\mathbf{W}} v^{-\sum\limits_{k<\ell}\langle\alpha_{j_\ell},\alpha_{j_k}\rangle}
|\mathcal{F}_\mathbf{j}(V)|\cdot\theta_i(\mathbf{j}).
\end{gather}
A~simple inspection using that $\theta_i(\mathbf{j})=0$ unless $\mathbf{j}$ ends with the letter~$i$ shows that
equations~\eqref{eq:derivation_RHS} and~\eqref{eq:derivation_LHS} agree.
The proof of (b) is similar.
\end{proof}

Our main conjecture of this note is that the quantum shuf\/f\/le character~$\Omega$ relates the representation theory of
$(Q,\mathbf{d})$ with the representation theory of KLR (quiver Hecke) algebras.
We will present necessary background and our conjectures in the next section.
First we draw a~useful conclusion from the existence of the homomorphism~$\Omega$.
\begin{Corollary}
\label{cor:counting polynomial}
For any rigid representation~$V$ of $(Q,\mathbf{d})$ and word $\mathbf{j}$, the variety $\mathcal{F}_\mathbf{j}(V)$ has
a~counting polynomial.
\end{Corollary}

\begin{proof}
Let~$V$ be a~rigid representation of $(Q,\mathbf{d})$ over a~f\/inite f\/ield $\mathbb{F}$ and write $\nu\in\mathcal{Q}$ for
the dimension vector of~$V$, a~positive root in~$\Phi$.
Since~$V$ is rigid, $[V]^*\in\mathcal{A}_v[N]\subset\mathcal{H}^*(Q,\mathbf{d})$, where $v=\sqrt{|\mathbb{F}|}$.
In fact, following~\cite{chen-xiao} and~\cite{jiang-sheng-xiao} there exists a~recursively computable element
$x_\nu\in\mathcal{A}_v[N]$, where~$v$ is a~formal parameter, which specializes to $[V]^*$ under the map $v\mapsto\sqrt{|\mathbb{F}|}$.
Since the recursion computing $x_\nu$ is described purely in terms of the root datum, for an arbitrary f\/inite f\/ield
$\mathbb{F}'$ the same element $x_\nu$ specializes under the map $v\mapsto\sqrt{|\mathbb{F}'|}$ to the class of the
unique indecomposable rigid representation of $(Q,\mathbf{d})$ over $\mathbb{F}'$ with dimension vector~$\nu$.

We may alternatively consider the image of $x_\nu$ under the embedding $\iota:\mathcal{A}_v[N]\to g^*$.
By the commutativity of the left hand triangle in~\eqref{dia:feigin_tetrahedron}, $\iota(x_\nu)$ specializes to
$\Omega([V]^*)$ under the map $v\mapsto\sqrt{|\mathbb{F}|}$.
Now recalling the def\/inition of $\Omega([V]^*)$ we see that (up to a~power of~$v$) the counting polynomial of
$\mathcal{F}_\mathbf{j}(V)$ is the coef\/f\/icient of $\mathbf{j}$ in $\iota(x_\nu)$.
As above, this polynomial is independent of the choice of specialization.
\end{proof}

\section{KLR algebras}\label{sec:KLR}

Recall the root lattice $\mathcal{Q}$ associated to the symmetrizable $I\times I$ Cartan matrix~$A$ with symmetri\-zers
$\mathbf{d}=(d_i: i\in I)$.
Write $\mathcal{Q}_+=\bigoplus\limits_{i\in I}\mathbb{Z}_{\ge0}\alpha_i$ for the positive cone in $\mathcal{Q}$
generated by the simple roots $\{\alpha_i\}_{i\in I}$.
The KLR (quiver Hecke) algebra $R=\bigoplus\limits_{\nu\in\mathcal{Q}_+} R_\nu$ is def\/ined diagrammatically~by
Khovanov--Lauda~\cite{khovanov-lauda1,khovanov-lauda2}
and symbolically in the independent work of
Rouquier~\cite{rouquier1,rouquier2}.
We follow the symbolic approach as presented in~\cite{kleshchev-ram2}.

Let $\Bbbk$ denote any f\/ield of characteristic zero.
The relations in $R_\nu$ will depend on the polynomials
\begin{gather*}
Q_{ij}(u,v):=
\begin{cases}
0 & \text{if $i=j$},
\\
1 & \text{if $a_{ij}=0$},  
\\
\varepsilon_{ij}(u^{-a_{ij}}-v^{-a_{ji}}) & \text{if $a_{ij}<0$},
\end{cases}
\end{gather*}
where, for simple representations $S_i$ and $S_j$ of $(Q,\mathbf{d})$, we set $\varepsilon_{ij}=1=-\varepsilon_{ji}$ if
$\operatorname{Ext}^1(S_i,S_j)\ne0$ and $\varepsilon_{ij}=0$ if
$\operatorname{Ext}^1(S_i,S_j)=\operatorname{Ext}^1(S_j,S_i)=0$.
Suppose the sum of the entries in~$\nu$ is~$d$.
Then the \emph{KLR $($quiver Hecke$)$ algebra} $R_\nu$ is the $\Bbbk$-algebra with polynomial generators $y_1,\ldots,y_d$,
braid-like generators $\psi_1,\ldots,\psi_{d-1}$, and orthogonal idempotents $e_\mathbf{j}$, where
$\mathbf{j}\in\mathbf{W}^\nu$ ranges over words with $|\mathbf{j}|=\nu$ (note that all such words have length~$d$).
These generators satisfy the following relations:
\begin{gather*}
\sum\limits_{\mathbf{j}\in\mathbf{W}^\nu}e_\mathbf{j}=1_\nu,
\qquad
y_ke_\mathbf{j}=e_\mathbf{j} y_k,
\qquad
y_k y_\ell=y_\ell y_k,
\\
\psi_k e_\mathbf{j}=e_{\sigma_k\mathbf{j}}\psi_k,
\qquad
(y_k\psi_\ell-\psi_\ell
y_{\sigma_\ell(k)})e_\mathbf{j}=\delta_{j_\ell,j_{\ell+1}}(\delta_{k,\ell+1}-\delta_{k,\ell})e_\mathbf{j},
\\
\psi_k^2 e_\mathbf{j}=Q_{j_k,j_{k+1}}(y_k,y_{k+1}) e_\mathbf{j},
\qquad
\psi_k\psi_\ell=\psi_\ell\psi_k \text{for $|k-\ell|>1$,}
\\
(\psi_{k+1}\psi_k\psi_{k+1}-\psi_k\psi_{k+1}\psi_k)
e_\mathbf{j}=\delta_{j_k,j_{k+2}}\frac{Q_{j_k,j_{k+1}}(y_{k+2},y_{k+1})-Q_{j_k,j_{k+1}}(y_k,y_{k+1})}{y_{k+2}-y_k}
e_\mathbf{j}.
\end{gather*}
The KLR algebra $R_\nu$ admits a~$\mathbb{Z}$-grading given~by
\begin{gather*}
\deg(e_\mathbf{j})=0,
\qquad
\deg(y_k e_\mathbf{j})=2d_{j_k},
\qquad
\deg(\psi_k e_\mathbf{j})=-d_{j_k}a_{j_k j_{k+1}}.
\end{gather*}

Consider the category $R-\operatorname{mod}$ of all f\/inite-dimensional graded~$R$-modules.
This category decomposes into the blocks $R_\nu-\operatorname{mod}$ ($\nu\in\mathcal{Q}_+$) consisting of
f\/inite-dimensional graded $R_\nu$-modules.

Using the natural embeddings $R_\nu\otimes R_\mu\hookrightarrow R_{\nu+\mu}$ we get induction and restriction functors
\begin{gather*}
\operatorname{Ind}_{\nu,\mu}^{\nu+\mu} : \ R_\nu\otimes R_\mu-\operatorname{mod}\to R_{\nu+\mu}-\operatorname{mod},
\\
\operatorname{Res}_{\nu,\mu}^{\nu+\mu} : \ R_{\nu+\mu}-\operatorname{mod}\to R_\nu\otimes R_\mu-\operatorname{mod}.
\end{gather*}
These are adjoint in the following sense:{\samepage
\begin{gather*}
\operatorname{Hom}_{R_{\nu+\mu}}(\operatorname{Ind}_{\nu,\mu}^{\nu+\mu}M,N)\cong\operatorname{Hom}_{R_\nu\otimes
R_\mu}(M,\operatorname{Res}_{\nu,\mu}^{\nu+\mu}N)
\end{gather*}
for every $M\in R_\nu\otimes R_\mu-\operatorname{mod}$ and $N\in R_{\nu+\mu}-\operatorname{mod}$.}

For $\nu,\mu\in\mathcal{Q}_+$ and $M\in R_\nu-\operatorname{mod}$, $N\in R_\mu-\operatorname{mod}$ we def\/ine $M\circ
N\in R_{\nu+\mu}-\operatorname{mod}$~by
\begin{gather*}
M\circ N=\operatorname{Ind}_{\nu,\mu}^{\nu+\mu}(M\boxtimes N).
\end{gather*}
With this product the Grothendieck group
$\mathcal{K}(R-\operatorname{mod})=\bigoplus_{\nu\in\mathcal{Q}_+}\mathcal{K}(R_\nu-\operatorname{mod})$ generated~by
isomorphism classes of f\/inite dimensional graded~$R$-modules has the structure of an associative unital ring. 
The identity element is given by the isomorphism class of $R_0$ thought of as an $R_0$-module and multiplication is
given by $[M][N]=[M\circ N]$ for $M\in R_\nu-\operatorname{mod}$ and $N\in R_\mu-\operatorname{mod}$.

For a~f\/inite dimensional graded vector space $V=\bigoplus_{n\in\mathbb{Z}} V_n$ write $\operatorname{dim}_v
V=\sum\limits_{n\in\mathbb{Z}}(\operatorname{dim} V_n)v^n$ for its \emph{graded character}.
The algebra~$R$ is graded and thus its Grothendieck ring admits an action of $\mathbb{Z}[v^{\pm1}]$, where~$v$ acts by grading shift.
There is a~$\mathbb{Z}[v^{\pm1}]$-linear character map $\operatorname{ch}_v$ from $\mathcal{K}(R-\operatorname{mod})$ to
the quantum shuf\/f\/le algebra $g^*$ given~by
\begin{gather*}
[M]\mapsto \operatorname{ch}_v(M):=\sum\limits_{\mathbf{j}\in\mathbf{W}} \operatorname{dim}_v M_\mathbf{j}\cdot\mathbf{j},
\end{gather*}
where $M_\mathbf{j}=e_\mathbf{j} M$ and the sum is f\/inite since~$M$ is a~f\/inite-dimensional representation.
The following theorem provides the monoidal categorif\/ication of quantum groups.
\begin{Theorem}[\cite{khovanov-lauda2,rouquier2,varagnolo-vasserot}]
The map $\operatorname{ch}_v: \mathcal{K}(R-\operatorname{mod})\to g^*$ is an injective algebra homomorphism, i.e.
\begin{gather*}
\operatorname{ch}_v(M\circ N)=\operatorname{ch}_v(M)\circ\operatorname{ch}_v(N),
\end{gather*}
and its image is precisely $\mathcal{A}_v[N]$.
Moreover, in symmetric types the irreducible characters correspond to elements of the dual canonical basis.
\end{Theorem}

Thus it is natural that a~great deal of ef\/fort has been put forth to describe and construct the irreducible
representations of KLR algebras.
Currently the main approach is to realize them as the irreducible heads of certain ``standard modules'' which admit
a~combinatorial description.
This approach has been pursued by Kleshchev--Ram~\cite{kleshchev-ram2}, McNamara~\cite{mcnamara}, and
Benkart--Kang--Oh--Park~\cite{benkart-kang-oh-park}
in f\/inite types and by Kleshchev~\cite{kleshchev} in the af\/f\/ine case.
An inductive approach to classifying irreducible~$R$-modules has been suggested by Lauda--Vazirani 
in~\cite{lauda-vazirani}.
We propose an entirely new construction of characters for certain irreducible~$R$-modules using the representation
theory of $(Q,\mathbf{d})$ and the quantum shuf\/f\/le character.

\begin{Conjecture}
\label{conj:rigid to irreducible}
For any rigid representation~$V$ of $(Q,\mathbf{d})$, the quantum shuffle character $\Omega([V]^*)$ produces the
character of an irreducible KLR representation.
\end{Conjecture}

The following example provides some small evidence for Conjecture~\ref{conj:rigid to irreducible} for arbitrary valued quivers $(Q,\mathbf{d})$.
\begin{Example}
It is shown in~\cite[Theorem 6.10]{lauda-vazirani} that for $\nu=-a_{ij}\alpha_i+\alpha_j$ ($i\ne j$) the block
$R_\nu-\operatorname{mod}$ contains a~unique irreducible module $L(i^{-a_{ij}},j)$ with
$\operatorname{ch}_v(L(i^{-a_{ij}},j))=[-a_{ij}]_i^!(i^{-a_{ij}},j)$ and a~unique irreducible module $L(j,i^{-a_{ij}})$
with $\operatorname{ch}_v(L(j,i^{-a_{ij}}))=[-a_{ij}]_i^!(j,i^{-a_{ij}})$.

Now consider the full subquiver $(Q_{ij},(d_i,d_j))$ of $(Q,\mathbf{d})$ on vertices~$i$ and~$j$ (without loss of
generality assume the arrows point from~$i$ to~$j$).
Any representation of $(Q_{ij},(d_i,d_j))$ is automatically a~representation of $(Q,\mathbf{d})$ and rigidity is carried
through this equivalence.
The quiver $(Q_{ij},(d_i,d_j))$ has a~unique non-simple injective representation $I_j$, the injective hull of the simple
$S_j$, and its dimension vector is precisely $-a_{ij}\alpha_i+\alpha_j$.
We aim to show that $\Omega([I_j]^*)=[-a_{ij}]_i^!(i^{-a_{ij}},j)$ is the character of the irreducible
$L(i^{-a_{ij}},j)$.
Indeed, notice that $I_j$ only admits f\/lags of type $(i^{-a_{ij}},j)$ and the coef\/f\/icient of $(i^{-a_{ij}},j)$ in
$\Omega([I_j]^*)$ is $v^{-{-a_{ij}\choose 2}}(-a_{ij})_i^!=[-a_{ij}]_i^!$ where we used that  
$\langle\alpha_j,\alpha_i\rangle=0$.

Similarly, replacing~$Q$ by $Q^{\rm op}$ and $Q_{ij}$ by $Q_{ij}^{\rm op}$ we see by a~similar argument that the unique
non-simple projective representation $P_j$, the projective cover of the simple $S_j$, satisf\/ies
$\Omega([P_j]^*)=[-a_{ij}]_i^!(j,i^{-a_{ij}})$ is the character of the irreducible $L(j,i^{-a_{ij}})$.
\end{Example}

For symmetric types, i.e.~equally valued quivers~$Q$, Conjecture~\ref{conj:rigid to irreducible} can be seen as
a~consequence of the geometric construction of Varagnolo--Vasserot and Webster~\cite{varagnolo-vasserot,webster}.

\begin{Theorem}
\label{th:shuffle character is KLR character}
Conjecture~{\rm \ref{conj:rigid to irreducible}} holds in the case $d_i=d_j$ for all $i,j\in Q_0$.
\end{Theorem}
\begin{proof}
We will begin by f\/ixing some notation roughly following~\cite{varagnolo-vasserot}.
Let~$V$ be a~rigid representation of~$Q$ and write $\nu=|V|$ for its dimension vector.
Let $\mathcal{E}_\nu$ be the space of all~$Q$-representations on the $\mathbb{N}^{Q_0}$ graded vector space
underlying~$V$.
The group $\mathcal{G}_V=\prod\limits_{i\in Q_0} GL(V_i)$ acts on $\mathcal{E}_\nu$ by the standard conjugation action.
Write $\mathcal{D}_{G_V}(\mathcal{E}_\nu)$ for the $\mathcal{G}_V$-equivariant
bounded derived category of complexes of sheaves on $\mathcal{E}_\nu$.
Let $\widetilde{\mathcal{F}}_\nu$ denote the variety of all pairs $(x,\phi)$ where $x\in\mathcal{E}_\nu$
is a~$Q$-representation and~$\phi$ is a~complete f\/lag of subrepresentations.
The group~$\mathcal{G}_V$ acts on $\widetilde{\mathcal{F}}_\nu$ and the projection
$\pi:\widetilde{\mathcal{F}}_\nu\to\mathcal{E}_\nu$ is a~$\mathcal{G}_V$-equivariant proper morphism.
The va\-riety~$\widetilde{\mathcal{F}}_\nu$ can be decomposed as a~disjoint union of isotypic components
$\widetilde{\mathcal{F}}_\nu=\bigsqcup\limits_{\mathbf{j}\in\mathbf{W}^\nu} \widetilde{\mathcal{F}}_\mathbf{j}$ consisting
of pairs $(x,\phi)$, where~$\phi$ is a~f\/lag of type $\mathbf{j}$; this decomposition respects the action of~$\mathcal{G}_V$.
Write $L_\mathbf{j}\in\mathcal{D}_{G_V}(\mathcal{E}_\nu)$ for the normalized pushforward under the map~$\pi$ of the
constant sheaf on~$\mathcal{F}_\mathbf{j}$, i.e.~shift the pushforward up by the dimension of~$\mathcal{F}_\mathbf{j}$.

\begin{Theorem}[\protect{\cite[Theorem 3.6]{varagnolo-vasserot}}]
The KLR-algebra $R_\nu$ is isomorphic to the Ext-algebra\linebreak
$\operatorname{Ext}_{\mathcal{G}_V}^\bullet(L_\nu,L_\nu)$,
where $L_\nu=\bigoplus\limits_{\mathbf{j}\in\mathbf{W}^\nu}L_\mathbf{j}$.
Under this isomorphism, the idempotent~$e_\mathbf{j}$ corresponds to the identity map on~$L_\mathbf{j}$.
\end{Theorem}

From a~$\mathcal{G}_V$-variety $\mathcal{Y}$ and a~$\mathcal{G}_V$-equivariant morphism
$\alpha:\mathcal{Y}\to\mathcal{E}_\nu$ we can construct the f\/iber product $\bigsqcup_\mathbf{j}
\mathcal{Y}\times_{\mathcal{E}_\nu}\mathcal{F}_\mathbf{j}$ whose $\mathcal{G}_V$-equivariant homology can be identif\/ied
with $\operatorname{Ext}_{\mathcal{G}_V}^\bullet(L,L_\nu)$, where~$L$ is the pushforward of the constant sheaf on
$\mathcal{Y}$ along the morphism~$f$ shifted up by the dimension of $\mathcal{Y}$, in particular there is an action of
$\operatorname{Ext}_{\mathcal{G}_V}^\bullet(L_\nu,L_\nu)$ on the homology.
Take $\mathcal{Y}=\mathcal{G}_V$ and let~$\alpha$ be the action map $\alpha(g)=g\cdot V$ whose image is identif\/ied with
the dense subset of $\mathcal{E}_\nu$ consisting of representations isomorphic to the original rigid representation~$V$.
In this case the f\/iber product above is identif\/ied with $\bigsqcup_\mathbf{j} \mathcal{F}_\mathbf{j}(V)$, where the
$e_\mathbf{j}$ weight space of the $R_\nu$ action is exactly the homology of $\mathcal{F}_\mathbf{j}(V)$.
By Corollary~\ref{cor:counting polynomial} and~\cite[Theorem~6.1.2]{hausel-rodriguez}, the Poincar\'e polynomial of
$\mathcal{F}_\mathbf{j}(V)$ is equal to its counting polynomial.
Taking into account our choice of shift and that the dimension of $\mathcal{F}_\mathbf{j}(V)$ is
$\sum\limits_{k<\ell}\langle\alpha_{j_\ell},\alpha_{j_k}\rangle$, it follows that $\Omega([V]^*)$ is the character of
the~$R_\nu$-module $H_*^{G_V}\big(\bigsqcup_\mathbf{j} \mathcal{F}_\mathbf{j}(V)\big)$.

It only remains to see that this module is irreducible.
First note that $\bigsqcup_\mathbf{j} \mathcal{F}_\mathbf{j}(V)$ is exactly the f\/iber of the map~$\pi$ over the point
$V\in\mathcal{E}_\nu$.
It follows (as in~\cite[equation~(8.6.24)]{chriss-ginzburg}) that the Jordan--H\"older multiplicity of a~simple $R_\nu$-module
in any composition series of $H_*^{G_V}\big(\bigsqcup_\mathbf{j} \mathcal{F}_\mathbf{j}(V)\big)$ is determined by the stalk over~$V$ of the
corresponding IC sheaf appearing in $L_\nu$.
But \mbox{$\mathcal{G}_V\cdot \{V\}\subset\mathcal{E}_\nu$} is a~locally closed, smooth, irreducible dense subvariety which is
$\mathcal{G}_V$-equivariantly simply connected.
It follows (see~\cite[Section~1.5]{lusztig2}) that there is a~unique perverse sheaf on $\mathcal{E}_\nu$
whose restriction to $\mathcal{G}_V\cdot\{V\}$ is nonzero.
In fact, the restriction will be the constant local system on $\mathcal{G}_V\cdot \{V\}$ shifted by the dimension of the
orbit.
Thus there is a~unique simple $R_\nu$-module appearing in $H_*^{G_V}\big(\bigsqcup_\mathbf{j}
\mathcal{F}_\mathbf{j}(V)\big)$ and it appears with multiplicity one, as desired.
\end{proof}

\subsection{Relation to dual canonical basis conjecture\\ for acyclic quantum cluster algebras}
The corollary below establishes the dual canonical basis conjecture for symmetric types re\-co\-ve\-ring a~result of
Kimura--Qin~\cite{kimura-qin}.
\begin{Corollary}
Assume~$A$ is a~symmetric Cartan matrix and let $w=c^2$ for a~Coxeter element $c\in W$.
Each non-initial cluster monomial of $\mathcal{A}_v[N^w]$ is contained in the dual canonical basis.
\end{Corollary}
\begin{proof}
Let $\mathbf{i}_0$ be a~reduced word for~$c$ and consider any valued quiver $(Q,\mathbf{d})$ which admits source adapted
sequence $\mathbf{i}_0$.
Write $\mathbf{i}=(\mathbf{i}_0,\mathbf{i}_0)$.
From Theorem~\ref{th:quantum cluster structure} we know that there is a~quantum cluster structure on
$\mathcal{A}_v[N^w]$ and by Theorem~\ref{th:qcc} the elements
$X_V=\widetilde\Psi_\mathbf{i}([V]^*)\in\widetilde\Psi_\mathbf{i}(\mathcal{A}_v[N^w])$ for~$V$ a~rigid representation of
$(Q,\mathbf{d})$ exactly correspond to non-initial quantum cluster monomials.
But we have the factorization $\widetilde\Psi_\mathbf{i}=\bar\Psi_\mathbf{i}\circ\Omega$ from Theorem~\ref{th:psi
factorization} which gives $X_V=\bar\Psi_\mathbf{i}(\Omega([V]^*))$ is the image of the dual canonical basis element $\Omega([V]^*)$.
\end{proof}

In the symmetrizable types it is known that dual canonical basis elements need not have positive multiplicative
structure constants, while this is still suspected to hold for skew-symmetri\-zable quantum cluster algebras.
The validity of Conjecture~\ref{conj:rigid to irreducible} for symmetrizable types would suggest the following
ref\/inement of the dual canonical basis conjecture for skew-symmetrizable quantum cluster algebras.
\begin{Conjecture}
For any Weyl group element $w\in W$ there is a~quantum cluster structure on~$\mathcal{A}_v[N^w]$, where cluster monomials
correspond to irreducible characters of KLR representations.
\end{Conjecture}

\section{Consequences and further conjectures}\label{sec:generalizations}

We conclude this note with some ideas of possible continuations or future directions for this approach.

\subsection{Conjectures on irreducible quiver Hecke modules\\ and quantum cluster algebras}
Fix a~Weyl group element $w\in W$.
It was shown by Lusztig in~\cite{lusztig, lusztig3} that the dual canonical basis respects the surjection
$\mathcal{A}_v[N]\onto\mathcal{A}_v[\overline{N^w}]$, i.e.~there is a~unique subset of the dual canonical basis which
induces a~basis in $\mathcal{A}_v[\overline{N^w}]$, we call this the dual canonical basis of
$\mathcal{A}_v[\overline{N^w}]$.
A~dual canonical basis element~$b$ is called \emph{$w$-adapted} if~$b$ corresponds to an element in the dual canonical
basis of $\mathcal{A}_v[\overline{N^w}]$.
\begin{Lemma}
A~dual canonical basis element~$b$ is~$w$-adapted if and only if $\Psi_\mathbf{i}(b)\ne0$ for a~reduced word
$\mathbf{i}=(i_1,\ldots,i_m)$ of~$w$.
\end{Lemma}
\begin{proof}
This is an immediate consequence of Theorem~\ref{th:ideal} and the discussion above.
\end{proof}

Since in symmetric types irreducible characters for the KLR algebra give dual canonical basis elements, we will say for
arbitrary symmetrizable types that an irreducible representation~$L$ of~$R$ is \emph{$w$-adapted} if
$\overline\Psi_\mathbf{i}(\operatorname{ch}_v(L))\ne0$ for a~reduced word $\mathbf{i}=(i_1,\ldots,i_m)$ of~$w$.
Following the cluster algebra philosophy we will call an irreducible character a~\emph{$w$-adapted coefficient} of its
image under $\overline\Psi_\mathbf{i}$ is a~single monomial and a~\emph{$w$-adapted non-coefficient} otherwise.
\begin{Conjecture}
Each~$w$-adapted coefficient is invertible in $\mathcal{A}_v[N^w]$ and the set of all~$w$-adapted non-coefficients forms
a~basis of $\mathcal{A}_v[N^w]$ over the Laurent ring of coefficients.
\end{Conjecture}

From now on we will restrict ourselves to~$w$-adapted irreducible representations.
To describe an analogue of clusters we make the following observation.
Cluster monomials have the property that the product of any two cluster monomials from the same cluster is (up to
a~power of~$v$) again a~cluster monomial.
Thus we make the following translation to the~$R$-module language.
\begin{Definition}
We def\/ine a~\emph{stably irreducible family} of~$R$-modules to be a~set $\mathbf{L}$ of irreducible modules such that
for all $L,L'\in\mathbf{L}$, $L\circ L'$ is irreducible and is (up to a~grading shift) an element of $\mathbf{L}$.
\end{Definition}
Kleshchev--Ram~\cite{kleshchev-ram2} have shown that arbitrary self-inductions of cuspidal irreducible representations
are again irreducible, in particular this shows that stably irreducible families exist.
Conjecture~\ref{conj:rigid to irreducible} also implies the existence of stably irreducible families, namely the
irreducible~$R$-modules corresponding to summands of any f\/ixed maximal rigid representation of $(Q,\mathbf{d})$.
Since our cluster monomials should be monomials in~$m$ variables we make the following def\/inition.

\begin{Definition}
For a~stably irreducible family $\mathbf{L}$ we call an element $L\in\mathbf{L}$ \emph{atomic} if there do not exist
$L',L''\in\mathbf{L}$ such that $L'\circ L''=L$.
Call a~stably irreducible family $\mathbf{L}$ \emph{$k$-generated} if there exists a~minimal collection of atomic
elements $\{L_1,\ldots,L_k\}\subset\mathbf{L}$ so that (up to a~grading shift) every element of $\mathbf{L}$ is of the
form $L_{j_1}\circ\dots\circ L_{j_r}$ ($1\le j_\ell\le k$) for some $r\ge1$.
In this case we will write $\mathbf{L}=\langle L_1,\ldots,L_k\rangle$.
\end{Definition}
Since the stably irreducible families are naturally ordered by inclusion, we may speak about maximal stably irreducible
families of~$w$-adapted representations.
\begin{Conjecture}
\label{conj:max-gen}
Every maximal stably irreducible family $($of~$w$-adapted representations$)$ is~$m$-generated $($where~$m$ is the length of~$w)$.
\end{Conjecture}

Thus we propose that maximal stably irreducible families should correspond to the clusters of $\mathcal{A}_v[N(w)]$,
where cluster variables are given by the atomic elements.

\begin{Remark}
The restriction to~$w$-adapted representations in Conjecture~\ref{conj:max-gen} seems to be essential.
Indeed, in general for inf\/inite-types by allowing arbitrary representations there should exist stably irreducible
families generated by any f\/inite number of atomic elements.
\end{Remark}

Having in mind the analogy with the hypothetical cluster structure on $\mathcal{A}_v[N^w]$ we make the following conjectures.
\begin{Conjecture}
\label{conj:coefficients}
Suppose each reduced expression for~$w$ contains every simple ref\/lection at least once.
Then there exist an~$n$-generated stably irreducible collection $\mathbf{C}=\langle C_1,\ldots,C_n\rangle$ such that $($up
to a~grading shift$)$ $\mathbf{C}\subset\mathbf{L}$ for every maximal stably irreducible family $\mathbf{L}$ and each
$C_i$ is atomic in $\mathbf{L}$.
\end{Conjecture}

As the notation suggests the $C_i$ should be thought of as coef\/f\/icients in the cluster associated to $\mathbf{L}$.
Assume Conjecture~\ref{conj:coefficients} holds.
\begin{Conjecture}
\label{conj:mutations}
Let $\mathbf{L}=\langle L_1,\ldots,L_{m-n},C_1,\ldots,C_n\rangle$ be a~maximal stably irreducible family and let
$L_k\in\mathbf{L}$ be an atomic non-coefficient.
Then there exists a~unique stably irreducible~$R$-module $L'_k\not\cong L_k$ such that $\mathbf{L}'=\langle
L_1,\ldots,L_{k-1},L'_k,L_{k+1},\ldots,L_{m-n},C_1,\ldots,C_n\rangle$ is a~maximal stably irreducible family different from $\mathbf{L}$.
\end{Conjecture}
A~proof of Conjecture~\ref{conj:mutations} would (essentially, modulo interpreting exchange matrices) complete the
monoidal categorif\/ication of the cluster structure on $\mathcal{A}_v[N^w]$ by KLR algebras.

\subsection*{Acknowledgements}
The author would like to thank Sasha Kleshchev for introducing him to KLR algebras and for leading him in this direction of research.
The author would also like to thank Arkady Berenstein for introducing him to Hall algebras and their beautiful properties.
Finally, special thanks need to be given to the anonymous referees for helping to solidify the proof of
Theorem~\ref{th:shuffle character is KLR character}.

\pdfbookmark[1]{References}{ref}
\LastPageEnding

\end{document}